\documentclass[11pt]{amsart}
\usepackage{amscd, amssymb, mathrsfs, mathabx, tikz-cd, amsthm, thmtools}
\usepackage{graphicx}

\usepackage{enumitem}

\author[Oh]{Jeongseok Oh}
\address{Department of Mathematics\\ 
Imperial College\\
London SW7 2AZ\\
United Kingdom}
\email{j.oh@imperial.ac.uk}

\usepackage{hyperref}

\input{xy}
\xyoption{all}

\newtheorem{Thm}{Theorem}[section]
\newtheorem{Conj}[Thm]{Conjecture}
\newtheorem{Prop}[Thm]{Proposition}
\newtheorem{Def}[Thm]{Definition}
\newtheorem{Def/Thm}[Thm]{Definition/Theorem}
\newtheorem{Cor}[Thm]{Corollary}
\newtheorem{Lemma}[Thm]{Lemma}

\theoremstyle{definition}
\newtheorem{Rmk}[Thm]{Remark}

\numberwithin{equation}{section}

\newcommand{\ot }{\otimes}
\newcommand{\ra }{\rightarrow}

\newcommand{\Hom }{{\mathrm{Hom}}}

\newcommand{\Spec}{{\mathrm{Spec}}}

\newcommand{\Res}{{\mathrm{Res}}}

\newcommand{\Sym}{{\mathrm{Sym}}}

\newcommand{\rank }{{\mathrm{rank}}}

\newcommand{\cO}{{\mathcal{O}}}
\newcommand{\cM}{{\mathcal{M}}}
\newcommand{\cL}{{\mathcal{L}}}
\newcommand{\cE}{{\mathcal{E}}}
\newcommand{\cF}{{\mathcal{F}}}
\newcommand{\cH}{{\mathcal{H}}}

\newcommand{\cC}{{\mathcal{C}}}

\newcommand{\cT}{{\mathcal{T}}}

\newcommand{\fC}{{\mathfrak{C}}}

\newcommand{\NN}{{\mathbb N}}

\newcommand{\PP }{{\mathbb P}}

\newcommand{\QQ }{{\mathbb Q}}
\newcommand{\CC }{{\mathbb C}}
\newcommand{\ZZ }{{\mathbb Z}}

\newcommand{\fS}{\mathfrak{S}}

\newcommand{\lan}{\langle}
\newcommand{\ran}{\rangle}

\newcommand{\cN}{\mathcal{N}}

\newcommand{\cS}{\mathcal{S}}

\newcommand\arXiv[1]{\href{http://arxiv.org/abs/#1}{arXiv:#1}}
\newcommand\mathAG[1]{\href{http://arxiv.org/abs/math/#1}{math.AG/#1}}

\begin{document}

\title{Quasimaps to GIT fiber bundles and Applications}

\begin{abstract}
In \cite{B}, Brown proved that the $I$-function of a toric fibration lies on the overruled Lagrangian cone of its $g=0$ Gromov--Witten theory, introduced by Coates and Givental \cite{CG}.
In this paper, we prove the theorem for partial flag variety fibrations.
To do so, we will construct new moduli spaces generalising the idea of Ciocan-Fontanine, Kim and Maulik \cite{CKM}.
\end{abstract}

\maketitle

\setcounter{tocdepth}{1}
\tableofcontents

\section{Introduction}
Gromov--Witten theory is a crucial part of one side of mirror symmetry. It can be encoded as a statement about a generating function of Gromov--Witten invariants.
As we may guess, a function generated by {\em genus zero ($g=0$ for short) invariants only} is relatively easier to study than the one with all invariants.
Interestingly it recovers higher genus invariants in some cases \cite{Gi} and may be of independent interest.

The {\em $J$-function} is the $1^{\mathrm{st}}$ order derivative of the $g=0$ generating function. We will see its precise definition and various other descriptions in the following section and Section \S \ref{PfChar}.
The {\em $I$-function} is known to be a counterpart of the $J$-function in mirror symmetry.
From mirror symmetry it follows that the $I$-function appears as a solution of differential equations which has a beautiful presentation as a hypergeometric series. This description can be used in several applications of the study of $g=0$ Gromov--Witten theory. However, in spite of the importance of the $I$-function it is difficult to compute the differential equations whose solution it is and further
it does not have a concrete definition.

In \cite{G2}, Givental described the $I$-function for a toric complete intersection as a generating function of invariants defined by integrations on certain moduli spaces.
Inspired by Givental's idea, Ciocan-Fontanine, Kim and Maulik constructed quasimap moduli spaces for a GIT quotient \cite{CKM} which suggested a concrete definition of the $I$-function in a greater generality \cite{CK0}. One advantage of this approach is that one can study the $I$-function {\em without} knowledge of the mirror.
The purpose of this article is to construct moduli spaces generalising the results to a fiber bundle with GIT quotient fibers.

Meanwhile, the relationship between $I$ and $J$-functions is highly non trivial even if we know what the $I$-function should be. That relationship is called {\em the mirror conjecture/theorem}.
Various mirror theorems have been proved since the seminal breakthrough by Givental in \cite{G2}.
In \cite{B}, Brown proved the mirror theorem for a fiber bundle whose fiber is a toric variety, thus proving a conjecture due to Elezi \cite{E}, see Conjecture \ref{Elezi} and Theorem \ref{Brown} below. 
Note that Elezi's conjecture was about projective bundles. Brown's result is in fact stronger than the original conjecture. His proof does not use moduli spaces to establish the conjecture. The main purpose of this article is to prove the mirror theorem for fiber bundles whose fiber is a flag variety using moduli spaces.

\subsection*{Gromov--Witten invariants}
Before explaining further, we would like to review Gromov--Witten invariants for a smooth projective variety $X$ over $\mathbb{C}$.
Roughly Gromov--Witten invariants count morphisms from curves to $X$ passing through given subvarieties of $X$. To be precise, we consider a collection of such morphisms
\begin{align} \label{modulii}
\mathcal{M}_{g,k}(X, \beta) :=  \{(C,f) \ | \ & C \text{ is a smooth curve of genus $g$ with $k$-marked}, \nonumber \\
\text{points, }&  f: C \rightarrow X \text{ is of degree $\beta$} \in H_2(X,\mathbb{Z}) \}/ \sim_{\text{isom}} 
\end{align} 
together with the evaluation maps
$$
\mathrm{ev}_a:\cM_{g,k}(X , \beta) \rightarrow X, \ \ \  (C, p_1, ..., p_k, f: C \rightarrow X) \mapsto f(p_a)
$$
at the marked points. Denoting by $X_a \in H^*(X, \mathbb{Q}[\![q]\!])$ the Poincar\'e dual of the $a^{\mathrm{th}}$ given subvariety, the invariant can be thought of as an intersection number 
$$
\deg\left( [\cM_{g,k}(X, \beta)] \ \cap \ \prod_{a=1}^k\  \text{ev}_a^*(X_a) \right) \ =\ \int_{[\cM_{g,k}(X, \beta)]} \ \prod_{a=1}^k\  \text{ev}_a^*(X_a) .
$$

A problem is that the space \eqref{modulii} may not be compact, thus the intersection number may not be defined.
Using stable maps, we can compactify the space \eqref{modulii}. The result is a Deligne--Mumford (DM) stack and it is usually denoted by $\overline{\mathcal{M}}_{g,k}(X, \beta)$. 
Another problem is that it may not be smooth, thus its fundamental class may not have the expected dimension.
However $\overline{\mathcal{M}}_{g,k}(X, \beta)$ is equipped with a natural perfect obstruction theory so that its {\em virtual} fundamental class $[\overline{\mathcal{M}}_{g,k}(X, \beta)]^{vir}$ is defined. Gromov--Witten invariant is then defined by
$$
\int_{[\overline{\mathcal{M}}_{g,k}(X, \beta)]^{vir}} \ \prod_{a=1}^k\  \text{ev}_a^*(X_a).
$$

\subsection*{$g=0$ Gromov--Witten theory} \label{g0GWth}
The $2^{\mathrm{nd}}$ homology classes $\beta$ defining nonempty $\overline{\mathcal{M}}_{g,k}(X, \beta)$ form a monoid in $H_2(X,\mathbb{Z})$. We denote it by $\mathrm{Eff}(X)$, or $\mathrm{Eff}$ for short.
One can define the group ring
$$
\QQ[\![q]\!]\ :=  \ \mathbb{Q}[\![\mathrm{Eff}]\!] \  =\ \lan \ q^\beta \ |\ \beta \in \mathrm{Eff}\ \ran,
$$
called {\em the Novikov ring}.
The psi class, $\psi_a$, is defined as the $1^{\mathrm{st}}$ Chern class of the line bundle on $\overline{\mathcal{M}}_{g,k}(X , \beta)$ formed by cotangent lines of $C$ at $p_a$.
Then {\em the $g=0$ descendant potential of $X$} is defined by
$$
\mathcal{F}_0\ :=\ \sum_{k=0}^\infty \sum_{\beta \in \text{Eff}} \frac{q^\beta}{k!} \int_{[\overline{\mathcal{M}}_{0,k}(X, \beta)]^{vir}} \prod_{a=1}^k \sum_{n=0}^\infty \text{ev}_a^*(t_n)\psi_a^n,
$$
where $t_0, t_1, ... \in H^*(X, \mathbb{Q}[\![q]\!])$ are formal variables.  
The function $\cF_0$ is one key ingredient of mirror symmetry.
For instance, mirror symmetry predicts the quantum cohomology ring of $X$ is isomorphic to the Jacobian ring of its mirror. Letting $t_0  :=  \sum_i t^i\gamma^i$ for a basis $\{\gamma^i\}$ of $H^*(X)$ with the dual basis $\{\gamma_i\}$, the quantum cohomology is defined  by $\left(H^*(X, \mathbb{Q}[\![q]\!]), \star_{t_0}\right)$ with the product
$$
\gamma^i \star_{t_0} \gamma^j \ = \ \partial_{t^i}\partial_{t^j}\partial_{t^k} \mathcal{F}_0|_{t_1=t_2=...=0} \cdot \gamma_k.
$$
Note that the restriction to $q=t_0=0$ of the quantum cohomology is the usual cohomology. 

The Lagrangian cone $\cL ag_X$ of $g=0$ Gromov--Witten theory of $X$ is roughly defined as the graph of the differential (in an infinite dimensional vector space)
$$
\cL ag_X\ :=\ \Gamma\left(d\cF_0 \right) \ \subset \ T^\vee_{H^*(X, \mathbb{Q})[z] \ot_\QQ \QQ[\![q]\!]},
$$
see \cite{CG,G3} for the precise definition.
It is a (origin-shifted) cone in the vector space.
Using the (non-trivial) identification of symplectic spaces
$$
T^\vee_{H^*(X, \mathbb{Q})[z] \ot_\QQ \QQ[\![q]\!]} \ \cong \ \left( H^*(X, \mathbb{Q})[z] \oplus \frac{1}{z}H^*(X, \mathbb{Q})[\![z^{-1}]\!] \right) \ot_\QQ \QQ[\![q]\!] =: \cH,
$$
its shifting $-z^{-1} \cL ag_X$ is spanned by the derivatives of ($z^{-1}$-expansion of) the $J$-function $J_X(-z, {\bf t})$ \cite[Proposition 1]{CG}, where
{\medmuskip=-2mu
\thinmuskip=-2mu
\thickmuskip=-2mu
\nulldelimiterspace=-1pt
\scriptspace=0pt
\begin{align} \label{JofF}
J_{X}(z, \ \textbf{t}) \ :=& \ 1\ +\ \frac{\textbf{t}}{z} \  + \sum_{ (k,\beta) \neq (0,0),(1,0) } \ \frac{q^\beta}{k!}  (\text{ev}_{k+1})_*\left(\frac{\prod_{a=1 \ }^{k}\text{ev}_a^*(\textbf{t})[\overline{\mathcal{M}}_{0,k+1}(X,\ \beta)]^{vir}}{z(z\ -\ \psi_{k+1})} \right). 
\end{align}}
\par \noindent

Computing $J_X$ is very difficult in general, but mirror symmetry predicts that it is related to a different, often more computable function known as the $I$-function. A modern form of the mirror conjecture/theorem (introduced by Givental) is to prove that an $I$-function lies on the Lagrangian cone.  This then allows one to recover the $J$-function via a complicated procedure of Birkhoff factorisation and change of variables.

\subsection*{Mirror theorem}
When $X$ is a GIT quotient $X = V /\!\!/ {\bf G}$, there is a nice description of the $I$-function $I_X$ \cite{CKM, CK0}. Ciocan-Fontanine, Kim and Maulik constructed another compactification $Q_{g,k}(X,\beta)$ of the space \eqref{modulii} by allowing morphisms $C \ra [V/{\bf G}]$ taking the generic point to $X$. These are called {\em quasimaps}. 
Then the $I$-function is defined by
$$
I_{X}(z,  \textbf{t}) \ :=\ \sum_{ (k,\beta) } \ \frac{q^\beta}{k!}  (\text{ev}_{k+1})_*\left(\frac{\prod_{a=1 \ }^{k}\text{ev}_a^*(\textbf{t})[Q_{0,k+1}(X, \beta)]^{vir}}{e_{\CC^*}(N_{k,\beta}^{vir})} \right), 
$$
where $N_{k,\beta}^{vir}$ denotes the virtual normal bundle of $Q_{0,k+1}(X,\ \beta)$ to a ``graph quasimap space'', see \eqref{graphJform} and Proposition \ref{Ishape} for graph spaces and virtual normal bundles. The input variable $z$ is the equivariant parameter of $\CC^*$.

When a torus acts on $X$, Ciocan-Fontanine and Kim proved the mirror theorem which is a generalisation of the result by Givental \cite{G2} for a complete intersection $X$ in a toric variety.

In \cite{B}, Brown defined the $I$-function for a fiber bundle whose fiber is a toric variety and proved the mirror theorem, see Theorem \ref{Brown}.

\subsection*{The main result} \label{MainSection}
The main goal of this article is to define the $I$-function and prove the mirror theorem for $X =\cF^{(n)}$, a fiber bundle over a smooth projective variety Y whose fiber is a partial flag variety.

We start with $F$ the total space of a sum of line bundles $\oplus_{j=1}^r L_j $ on $Y$ and the $n$-tuple of positive integers
$$
(r_1, ..., r_n) \ \in \ \ZZ^n, \ \  0< r_1 < ...<r_n<r_{n+1}:=r.
$$
Then we let $\pi: \mathcal{F}^{(n)} \rightarrow Y$ be a fiber bundle over $Y$ whose fiber at $y \in Y$ is the space of all collections of subspaces  
\begin{align*}
F_1 \subset F_2 \subset \cdots \subset F_n \subset F|_y, \ \ \text{rank}\; F_i=r_i.
\end{align*}
Our plan is to write the $I$-function in terms of $J$-function of $Y$ and the $I$-function of the fiber. So let
$$
\QQ[\![Q]\!] \ := \ \QQ[\![\mathrm{Eff}(Y)]\!]
$$
be the Novikov ring of $Y$ and $J_{D}$ be the $Q^{D}$-coefficient of $J_Y$. With a formal variable $u \in H^*(Y, \mathbb{Q})$, we can write 
$$J_Y (z, u)\ =\ \sum_{D} Q^{D}  J_{D}(z ,u).$$
Let $\mathcal{F}_1 \subset \cdots \subset \mathcal{F}_n \subset \mathcal{F}_{n+1}:= \pi^*F$ be the tautological bundles on $\mathcal{F}^{(n)}$.
For an effective class $D \in \mathrm{Eff}(Y)$, let $\ZZ^n_D$ be the collection of all $\beta \in \mathrm{Eff}(\cF^{(n)})$ satisfying $\pi_*\beta=D$.
Then $\ZZ^n_D$ is a subset of $\ZZ^n$ via $\beta \mapsto (\beta(\det \cF^{\vee}_i))_i \in \ZZ^n$.
Furthermore, $\cup_{D \in \mathrm{Eff}(Y)} \left( \{D\} \times \ZZ^n_D \right)$ forms a monoid in $\mathrm{Eff}(Y) \times \ZZ^n$ since so is $\mathrm{Eff}(\cF^{(n)})$. Let $\QQ[\![q,Q]\!]$ be the group ring defined by the monoid. It is the Novikov ring of $\cF^{(n)}$ and $\QQ[\![q,Q]\!] \ra \QQ[\![Q]\!]$ is induced by the projection of the monoid to $\mathrm{Eff}(Y)$.
Let $H_{i,l}$ be Chern roots of $\mathcal{F}_i^\vee$, $i=1, ..., n$, $l=1, ..., r_i$ and $H_{n+1,j}:= -\pi^*(c_1(L_j))$.
Then for formal variables $\textbf{t} = \sum_{i=1}^n t_i c_1(\cF_i ^\vee) \in H^{ 2}(\mathcal{F}^{(n)}, \mathbb{Q} )$ and $u \in H^*(Y, \mathbb{Q})$, we define the $I$-function for $\cF^{(n)}$
\begin{align}
\setlength\arraycolsep{0.05cm}
\label{flag}
I_{\cF^{(n)}} (z, {\bf{t}}, u)
:= \
&  e^{\frac{\textbf{t}}{z}} \sum_{D \in \text{Eff}(Y), d= (d_i) \in \ZZ^n_D} q^d Q^D e^{\sum_i t_i d_i} \pi^*(J_{D}(z, u)) \\
&  \sum_{\sum_l d_i^l =d_i} \prod_{i=1}^n  \left( \prod_{1 \leq l \neq l' \leq r_i}\frac{\prod_{s=-\infty}^{d_i^l -d_i^{l'}}(H_{i,l}-H_{i,l'}+sz) }{\prod_{s=-\infty}^{0}(H_{i,l}-H_{i,l'}+sz)}  \right.  \nonumber \\
 & \left. \prod_{\substack{1 \leq l \leq r_i, \\ 1 \leq l' \leq r_{i+1}}}\frac{\prod_{s=-\infty}^{0}(H_{i,l}-H_{i+1,l'}+sz)}{\prod_{s=-\infty}^{d_i^l -d_{i+1}^{l'}}(H_{i,l}-H_{i+1,l'}+sz)} \right), \nonumber
\end{align} 
which is an element in $\cH$, where $d^j_{n+1}:= - c_1(L_j) \cap D = - D(L_j)$. The summation over $\sum_i d^l_i = d_i$ is taken as follows.
Note that in Section \S \ref{rep:F}, we will express $\cF^{(n)}$ as a quotient $E /\!\!/ {\bf G}$ of a vector bundle $E$ on $Y$ by a linearly reductive group $\bf G$.
$(d^l_i)$ in \eqref{flag} is an effective class of its abelianisation $E /\!\!/ {\bf G}_T$, where ${\bf G}_T$ is a maximal torus of ${\bf G}$, whose image under the push forward of the stack morphism $[E/{\bf G}_T] \to [E/ {\bf G}]$ is $d_i$.
Hence the summation is finite.

\smallskip

Here is our main result.

\begin{Thm}
\label{Mthm}
The $I$-function $I_{\cF^{(n)}}(-z, {\bf{t}}, u)$ lies on $-z^{-1}\mathcal{L}ag_{\mathcal{F}^{(n)}}$.
\end{Thm}

\subsection*{Equivariant theory} \label{Sect:EqThy}
We prove Theorem \ref{Mthm} using a natural fiberwise action of a torus $\textbf{S} := (\CC^*)^r$ on $\cF^{(n)}$.
Let $\lambda$ be the equivariant parameters of $\bf S$.
The action defines the $\bf{S}$-equivariant $I$-function for $\cF^{(n)}$, denoted by $I^{\textbf{S}}_{\cF^{(n)}} (z, {\bf{t}},u, \lambda)$. 
Theorem \ref{Mthm} follows from its equivariant version by taking $\lambda \ra 0$. 

\begin{Thm} \label{Eq:Main}
The $I$-function $I^{\bf{S}}_{\cF^{(n)}}( -z, {\bf{t}}, u, \lambda)$ lies on $-z^{-1}\cL ag^{\bf{S}}_{\cF^{(n)}}$.
\end{Thm}

We prove Theorem \ref{Eq:Main} through a characterisation of $\cL ag ^{\bf{S}} _{\cF^{(n)}}$ (Theorem \ref{Char}).
The characterisation determines whether or not a function
$$
G \ \in \ \cH \; [\![{\bf t}, u]\!] \ \subset \ H^*_{\bf S} (\cF^{(n)}, \QQ)(z) \ot_\QQ \QQ[\![q, Q, {\bf t}, u]\!]
$$
lies on $-z^{-1}\cL ag ^{\bf{S}} _{\cF^{(n)}}[\![{\bf t}, u]\!]$. Here, for a rational function in $z$ we consider its $z^{-1}$-expansion when we need to view it as a series.

\begin{Thm}
\label{Char}
Suppose that $G$ satisfies the following properties.
\begin{enumerate} 
\item Initial condition in Section \S \ref{subsub:initial},
\item 
Recursion relation in Section \S \ref{subsub:rec},
\item 
Polynomiality condition in Section \S \ref{subsub:poly}.
\end{enumerate}
Then $G(-z, {\bf t}, u, \lambda)$ lies on  $-z^{-1}\cL ag^{\bf{S}} _{\cF^{(n)}}$ at each ${\bf t}$ and $u$.
\end{Thm}

\subsection*{Plan of the paper}
We explain properties in Section \S \ref{Sect:conditions} and prove Theorem \ref{Char} in Section \S \ref{ProofChar}. A brief strategy is as follows. 
We first produce a function, say $F$, associated with $G$ which lies on the Lagrangian cone. 
Then we show $G=F$ using the characterisation properties for $G$. 
In Section \S \ref{Sect:reverse} we discuss a converse of Theorem \ref{Char}.

We will show that $I^{\bf S}_{\cF^{(n)}}$ satisfies the characterisation properties in
Section \S  \ref{RecurSion} and \S  \ref{PolyI}.
In Section \S \ref{RecurSion}, we will compute the residues at {\em some} simple poles to check \eqref{Recursion:Imu} holds for the $I$-function. 
The quasimap spaces are introduced and studied in Section \S \ref{QuasiMapModuli}, which is quite independent from other sections. 
Then we write $I^{\bf{S}}_{\cF^{(n)}}$ in terms of the quasimap spaces and check the properties for it using the geometric interpretation in Section \S \ref{PolyI}.

As an application of our main result in Section \S \ref{Sect:App} we show how the $g = 0$ Gromov--Witten invariants of the base and total space are related to each other when $\cF^{(n)}$ is assumed to be Fano or Calabi--Yau.

\subsection*{Relationship to previous works: what is new and what is not}

When $Y = \text{Spec}\; \mathbb{C}$, Theorem \ref{Mthm} specialises the result by Bertram, Ciocan-Fontanine, and Kim \cite{BCK1, BCK2}.
When $n=r_1=1$, Theorem \ref{Mthm} proves Elezi's conjecture \cite{E} as a special case. 
\begin{Conj}[Elezi, \cite{E}] \label{Elezi}
Let $L_1 := \cO_Y$ and $L_j$, $j \neq 1$, be nef line bundles on $Y$ such that $-K_Y -\sum_j c_1(L_j)$ is ample.
Then we have 
$$
I_{\PP(\oplus_{j=1}^r L_j)}\ =\ J_{\PP(\oplus_{j=1}^r L_j)}.
$$
\end{Conj}
In \cite{B}, Brown proves Elezi's conjecture. He considers a toric fibration $\pi: \mathcal{T}:=\oplus_{j=1}^r L_j /\!\!/ (\mathbb{C}^*)^n \ra Y$, $n \leq r$. It has $r$-many toric divisors
$$
D_j \ :=\ \oplus_{i \neq j} L_i /\!\!/ (\mathbb{C}^*)^n \ \subset\ \mathcal{T}
$$
and $n$-many line bundles $P_i$ corresponding to the $i$-th factor of the torus $(\CC^*)^n$. For ${\bf t}=\sum_{i=1}^n t_i c_1(P_i)$, he defined the $I$-function
$$
\setlength\arraycolsep{0.05cm}
\begin{array}{c}
 I_{\cT} (z, {\bf t}, u)
\ := \    e^{\frac{\bf t}{z}} \displaystyle\sum_{\beta \in \mathrm{Eff}} q^\beta  e^{\int_\beta {\bf t}} \pi^*(J_{\pi_*(\beta)}(z, u)) \prod_{j=1}^r \frac{\prod_{s=-\infty}^{0}(D_j +sz)}{\prod_{s=-\infty}^{D_j.\beta}(D_j +sz)} . \\
\end{array}
$$

\begin{Thm}[Brown, \cite{B}] \label{Brown}
The $I$-function $I_{\cT}(-z, {\bf t}, u)$ lies on $-z^{-1}\mathcal{L}ag_\mathcal{T}$.
\end{Thm}

For Theorem \ref{Brown}, Brown also used a characterisation \cite[Theorem 2]{B}.
His characterisation has a fairly different aspect from ours because he used an asymptotic analysis to check if the $I$-function satisfies the characterisation properties.
Whereas, we will check it using the geometric interpretation of the $I$-function.

\smallskip

Our characterisation is motivated by \cite[Proposition 4.5]{G2} and \cite[Lemma 7.7.1]{CK0}.
Essentially, both \cite[Proposition 4.5]{G2} and \cite[Lemma 7.7.1]{CK0} assert that a function satisfying the recursion relation \S \ref{subsub:rec} and the polynomiality condition \S \ref{subsub:poly} lies on the Lagrangian cone. 

In \cite{G2, CK0}, the $I$-functions were interpreted in terms of suitable moduli spaces. The recursion and the polynomiality were checked by virtual localisation \cite{GP}.

\subsubsection*{Characterisation properties}

Since our target space is a fiber bundle over $Y \neq \Spec \; \CC$, the recursion and the polynomiality used in \cite{G2,CK0} for $Y = \Spec \; \CC$ are not enough. So they need to be modified. 
Since the recursion relation is a recursive argument with respect to the fiber direction, we need an initial condition of that recursive argument in terms of the base space $Y$.
This is the reason why the initial condition \S \ref{subsub:initial} is imposed on the characterisation theorem.
The polynomiality condition should be extended to any directional derivative along a vector at the origin on the cohomology of $Y$, see \eqref{eq:poly} for the precise statement.
These modifications are \emph{newly} designed.

\subsubsection*{Quasimap moduli spaces}
 
For a construction of moduli spaces, we have to keep two things in mind.
First, the invariants of such moduli spaces produce $I^{\bf S}_{\cF^{(n)}}$ as a generating function.
Second, the construction must be natural enough to check the characterisation properties for $I^{\bf S}_{\cF^{(n)}}$ easily.

An immediate approach of the construction may be able to be a direct generalisation of \cite{G2,CK0} --
construct a space of maps which project prestable maps to the base $Y$ and quasimaps \cite{CKM} to the fiber.
\emph{However, this approach is not good enough to meet the first purpose.}
So we will introduce a new idea which considers a space of maps which project prestable maps to $Y$ and quasimaps \emph{from contracted domain curves} to the fiber. 
In Section \S \ref{QuasiMapModuli}, we formalise this construction.

\begin{Rmk}
In Section \S \ref{RecurSion}, we will compute residues in the recursion relation without using moduli spaces.
But it can be computed by virtual localisation as well.
\end{Rmk}

\medskip

\subsection*{Acknowledgments}
I would like to thank Bumsig Kim for initially suggesting the problem, giving countless invaluable suggestions, support and advice.
I could not have finished the project without his help.

I thank Alexander Givental for teaching me $g=0$ mirror symmetry in terms of the formalism, kindly answering my questions, and precious advice.
I also thank Ionu\c t Ciocan-Fontanine for valuable discussions, suggestions and comments on an earlier draft and Bhamidi Sreedhar for suggestions and comments on the final draft.
Special thanks are due to the anonymous referees for their patience and inspiring comments, advice as well as Michel van Garrel, Daniel Kaplan, Mark Shoemaker, Pablo Solis and Jonathan Wise.

J. O. is partially supported by KIAS individual grant MG063002.

\section{Characterisation theorem}

In this section, we would like to explain characterisation properties and prove Theorem \ref{Char}. Since it characterises elements in the Lagrangian cone which is generated by the $J$-function, we start with a study of it. 

Indeed, though the characterisation can be generalised to any fiber bundle with a nice fiberwise torus action, we perfer to focus on $\cF^{(n)}$.

\subsection{$J$-function} \label{PfChar}
We have seen the definition of the $J$-function in \eqref{JofF}. Its equivariant version $J^{\bf{S}}_{\cF^{(n)}}$ is described in the same way.
Meanwhile there is an alternative expression of $J^{\bf{S}}_{\cF^{(n)}}$ providing a geometric meaning for the variable $z$.
Consider the so-called graph space 
\begin{align} \label{Gsp}
\overline{\cM}G_{0,k}(\cF^{(n)}, \beta)\ :=\  \overline{\mathcal{M}}_{0,k}\left(\cF^{(n)} \times \mathbb{P}^1, (\beta,1)\right).
\end{align}
Then the $\mathbb{C}^*$-action on $\mathbb{P}^1$, 
$$
\CC^* \times \PP^1 \ra \PP^1, \ \ \ \left(t,[x;y]\right) \mapsto [tx;y],
$$ 
defines an action on $\overline{\mathcal{M}}G_{0,k}(\cF^{(n)} , \beta)$.
If we set the convention $0=[0;1]$, $\infty=[1;0]$, then we obtain 
$$
e_{\CC^*}(T_0 \PP^1)=z ,\ \ \ e_{\CC^*}(T_\infty \PP^1)=-z,
$$ 
where $z=e_{\CC^*}(\CC_1)$ is the Euler class of the one dimensional $\CC^*$-representation with weight 1. 
On the domain curve for an element in $\overline{\cM}G_{0,k}(\cF^{(n)}, \beta)$, there is the \emph{unique} rational component $\PP^1$ which identically maps to the $\PP^1$-factor on the target space.
We call this component on the domain curve the {\em distinguished $\PP^1$}. 
By abuse of notation, we will denote this component simply by $\PP^1$ when there is no danger of confusion.
Let 
$F_{k,\beta} \subset \overline{\mathcal{M}}G_{0,k}(\cF^{(n)} , \beta)^{\mathbb{C}^*}$
be a component of the $\CC^*$-fixed locus, where $\infty \in \PP^1$ is neither a marked point nor a node.
There is an isomorphism $F_{k,\beta} \cong \overline{\cM}_{0,k+1}(\cF^{(n)}, \beta)$ unless $(k,\beta) = (0,0), (1,0)$,
defined by contracting $\PP^1$ to the extra marked point. We call it the {\em distinguished marked point}.
Since $\PP^1$ maps constantly to $\cF^{(n)}$, we can define an evaluation map
$\text{ev}_\bullet: F_{k,\beta} \ra \cF^{(n)}$
by the constant image of $\PP^1$.
It is $\text{ev}_\bullet = \text{ev}_{k+1}$ through the isomorphism $F_{k,\beta} \cong \overline{\cM}_{0,k+1}(\cF^{(n)}, \beta)$ if $(k,\beta) \neq (0,0), (1,0)$.
Hence $J^{\bf{S}}_{\cF^{(n)}}$ \eqref{JofF} can be rewritten as 
\begin{align} \label{graphJform}
\setlength\arraycolsep{0.05cm}
J^{\bf{S}}_{\cF^{(n)}}(\textbf{t})\ =\ \displaystyle\sum_{ k,\beta} \frac{q^\beta}{k!} (\text{ev}_{\bullet})_*\left(\frac{\prod_{a=1}^{k}\text{ev}_a^*(\textbf{t})[F_{k,\beta}]^{vir}}{e_{\CC^* \times {\bf{S}}}(N^{vir}_{F_{k,\beta}/  \overline{\mathcal{M}}G_{0,k}(\cF^{(n)}, \beta)})} \right) ,
\end{align}
where $N^{vir}$ denotes the virtual normal bundle with respect to the $\CC^*$-action.
$q^{\beta} \in \QQ[\![\mathrm{Eff}(\cF^{(n)})]\!]$ can be thought of as $q^{(\beta(\det \cF^\vee_i))_i}Q^{\pi_*\beta} \in \QQ[\![q,Q]\!]$ under the isomorphism $\QQ[\![\mathrm{Eff}(\cF^{(n)})]\!]  \cong  \QQ[\![q,Q]\!]$.
Note that for $(k,\beta) \neq (0,0), (1,0)$, we have
$$e_{\CC^* \times {\bf{S}}}(N^{vir}_{F_{k,\beta}/  \overline{\mathcal{M}}G_{0,k}(\cF^{(n)}, \beta)})\ =\ z(z-\psi_\bullet)$$
from the contributions of moving and smoothing the nodal point $\bullet=0 \in \PP^1$.

There is one more useful description of the $J$-function.
We define an operator 
\begin{eqnarray*}
& S^*_{\bf t}(z): \; H^*_{\bf{S}}({\cF^{(n)}}, \QQ)[z] \ot_\QQ \QQ[\![ q]\!]\ \rightarrow \ H^*_{\bf{S}}({\cF^{(n)}}, \QQ)(z) \ot_\QQ \QQ[\![ q]\!], \\
&
\ \gamma \; \mapsto\;  \gamma  + \sum_{(m,\beta) \neq (0,0)} \frac{q^\beta}{m!}  (\text{ev}_1)_* \left( \frac{\text{ev}_2^*(\gamma)\prod_{i=1}^m \text{ev}_{i+2}^*({\bf{t}})[\overline{\cM}_{0,2+m}({\cF^{(n)}}, \beta)]^{vir}}{z-\psi_1}\right).
\end{eqnarray*}
Letting $\{\gamma_i\}$ be a basis for $H^*_{\bf{S}}({\cF^{(n)}}, \QQ )$ and $t_i$ be a formal variable corresponding to $\gamma_i$,
one can check $S^*_{\bf{t}}(z)(\gamma_i) = z\partial_{t_i} J^{\bf{S}}_{\cF^{(n)}}$ using \eqref{JofF}.
On the other hand, by \cite[Proposition 1(i)]{CG}, we have
$$T_f \cL\ \cap\ \cL ag_{\cF^{(n)}}\ =\ zT_f \cL$$
at any point $f \in \cL ag_{\cF^{(n)}}$, where $T_f \cL$ denotes the tangent space to $\cL ag_{\cF^{(n)}}$ at $f$. Hence any directional derivative $z\partial f$ lies on  $zT_f \cL \subset \mathcal{L}ag_{\cF^{(n)}}$.
Applying it to $f=-zJ^{\bf S}_{\cF^{(n)}}(-z)$, one can check that $S^*_{\bf{t}}(-z)(\gamma)$ lies on $-z^{-1}\cL ag^{\bf S}_{\cF^{(n)}}$ for any $\gamma$.
For the special case $\gamma=1$, we have
\begin{align}\label{aaa}
S^*_{\bf{t}}(z)(1)\ =\ J^{\bf{S}}_{\cF^{(n)}}
\end{align}
by the string equation of Gromov--Witten theory \cite[(SE)]{G3}. 
Conversely, $\{S^*_{\bf{t}}(-z)(\gamma_i)\}_i$ spans $T_{-zJ^{\bf S}_{\cF^{(n)}}(-z)}\cL$ as a free $\QQ[\![q]\!][z]$-module \cite[Proposition 1]{CG}. This property is called {\em the reconstruction of the Lagrangian cone}.

\subsection{Characterisation properties} \label{Sect:conditions}
Now we explain characterisation properties for a function
$$
G \ \in \ \cH \; [\![{\bf t}, u]\!] \ \subset \ H^*_{\bf S} (\cF^{(n)}, \QQ)(z) \ot_\QQ \QQ[\![q, Q, {\bf t}, u]\!].
$$ 
Note that the $J$-function \eqref{JofF}, \eqref{graphJform}, \eqref{aaa} is defined on the whole cohomology ${\bf t} \in H^*(\cF^{(n)}, \QQ)$ whereas $G$ is the function on the restriction $H^2(\cF^{(n)}, \QQ) \oplus H^*(Y, \QQ) \subset H^*(\cF^{(n)}, \QQ)$. By abuse of notation, we use $\bf t$ for the variable on $H^2(\cF^{(n)}, \QQ)$. $u$ denotes the variable on $H^*(Y, \QQ)$.

\subsubsection{Initial condition} \label{subsub:initial}

Let $\mu: Y \hookrightarrow  \cF^{(n)}$ be the inclusion of an $\bf{S}$-fixed locus. We denote its image by $Y^\mu:= \mu(Y)$.
In \cite{CG}, the $N_{Y^\mu / \cF^{(n)}}$-twisted Gromov--Witten theory on $Y \cong Y^\mu$ is defined using the pairing
$$
\int_Y e^{-1}_{\bf S} \left(N_{Y^\mu / \cF^{(n)}} \right) \ a \cup b , \ \ \ \ \ a, \ b \ \in H^*_{\bf S}(Y) \ \cong \ H^*(Y)[\lambda].
$$
It gives rise to the twisted Lagrangian cone $\cL^\mu_Y$, see \cite[Section 7--10]{CG}.

Denoting by $G^\mu :=  \sum_{(d,D)} q^dQ^D \ G^\mu_{(d,D)}$ the pullback $\mu^*G$, we define
$$
G^\mu_Y \ :=\  \sum_{D} Q^D \ G^\mu_{(d_{D,\mu},D)} , \ \ d_{D, \mu}  := (D(\det \mu^*\cF_i^\vee))_i \ \in\ \ZZ^n.
$$
{\em The initial condition} for $G$ asserts that 
\begin{itemize}
\item $G|_{q=Q=0}=e^{({\bf t} + \pi^*u)/z}$, and 
\item $G^\mu_Y(-z, {\bf t}, u)$ lies on $-z^{-1}\cL^\mu_Y$. 
\end{itemize}

\subsubsection{Recursion relation} \label{subsub:rec}

Let
$\mu, \nu: Y \hookrightarrow \cF^{(n)}$
be $\bf{S}$-fixed loci in $\cF^{(n)}$ such that there exists a one-dimensional orbit connecting $Y^\mu$ and $Y^\nu$.
Then for each point in $Y^\mu$, there is a unique one-dimensional orbit connecting this point to a point in $Y^\nu$. 
The tangent space to a one-dimensional orbit at each point in $Y^\mu$ forms an $\bf S$-equivariant line bundle on $Y \cong Y^\mu$.
Let 
\begin{align*} 
\chi_{\mu,\nu}\ \in\ H^2_{\bf S}(Y , \QQ), \ \ \ d_{\mu,\nu}\ \in\ \mathrm{Eff}(\cF^{(n)})
\end{align*}
be its $1^{\mathrm{st}}$ Chern class and the class of a representative of the one-dimensional orbit, respectively.

For each integer $k>0$, since $Y^\mu \cong Y$, we have an embedding $Y \hookrightarrow \overline{\cM}_{0,2}(\cF^{(n)}, kd_{\mu,\nu})^{\bf{S}}$ defined by the $k$-coverings of orbits.
Note that the two marked points are corresponding to $0, \infty \in \PP^1$ the poles of orbits.
Let 
\begin{align*} 
N^{vir}_{\mu,\nu,k}\ :=\ N^{vir}_{Y/ \overline{\cM}_{0,2}(\cF^{(n)}, kd_{\mu,\nu})}
\end{align*}
be the virtual normal bundle to $Y$ in $\overline{\cM}_{0,2}(\cF^{(n)}, kd_{\mu,\nu})$.

We say that $G$ satisfies {\em the recursion relation} if each coefficient of $G^\mu$, as a series in $q$, $Q$, $\bf{t}$, and $u$, has
\begin{itemize}
\item finite order poles at $z=0$ and $z=\infty$, and
\item simple poles at $z=-\chi_{\mu,\nu}/k$ for $k \in \NN$, with residues given by
\begin{align} \label{Recursion:Imu}
\Res_{z=-\frac{\chi_{\mu,\nu}}{k}} G^\mu(z)dkz\ =\  \frac{q^{kd_{\mu,\nu}}}{e_{\bf{S}}(N^{vir}_{\mu,\nu,k})} e_{\bf{S}}(N_{Y^\mu / \cF^{(n)}}) G^\nu (-\frac{\chi_{\mu, \nu}}{k}).
\end{align}
\end{itemize}

\subsubsection{Polynomiality condition} \label{subsub:poly}
Let $\{\delta_j\}$ be a basis for $H^*(Y, \QQ)$ and $u_j$ be the corresponding variable in $u= \sum_j u_j \delta_j \in H^*(Y, \QQ)$. Also let $E_{i,\mu} :=  \det \cF^\vee_{i} \otimes \pi^*\mu^* \det \cF_{i} $ be the bundle on $\cF^{(n)}$. Then {\em the polynomiality condition} for $G$ says that 
\begin{itemize}
\item for each $\mu$ and $j$, the rational function
\begin{align} \label{eq:poly}
(z\partial_{u_j}  G^\mu(z,q) ,  G^\mu(-z, qe^{-z\sum_i y_iE_{i,\mu}} ) )_Y
\end{align}
has no poles in $z =0$, where $y_i$'s are formal variables. 
\end{itemize}

\subsection{Converse of Theorem \ref{Char}} \label{Sect:reverse}
The $J$-function $J^{\bf S}_{\cF^{(n)}}$ satisfies the initial condition (by \cite[Theorem 2]{B} and the reconstruction of $\cL^\mu_Y$) and the recursion relation (see \cite{G2,CK0, B}).

\begin{Prop} \label{Prop:polyJ}
$J^{\bf S}_{\cF^{(n)}}$ satisfies the polynomiality condition as well.
\end{Prop}

\begin{proof}
Recall from \eqref{graphJform} that $J^{\bf S}_{\cF^{(n)}}$ is a summation of the integrations over $F_{k,\beta}$'s. $F_{k,\beta}$ is a component in $\overline{\cM}G_{0,k+1}(\cF^{(n)}, \beta)^{\CC^*}$.
We consider other components 
$$
F_{k_1,\beta_1}^{k_2,\beta_2} \ \subset\ \overline{\cM}G_{0,k+1}({\cF^{(n)}}, \beta)^{\CC^*}, \ \ k_1 +k_2=k, \ \beta_1 + \beta_2 = \beta,
$$ 
where $k_1$-marked points, degree $\beta_1$ are concentrated on $0=[1;0] \in \mathbb{P}^1$ and $k_2$-marked points, degree $\beta_2$ are concentrated on $\infty \in \mathbb{P}^1$, as well. Note that \begin{align} \label{hahaiso}
F_{k_1,\beta_1}^{k_2,\beta_2} \ \cong F_{k_1,\beta_1} \times_{\cF^{(n)}} F_{k_2, \beta_2}.
\end{align}
Let $\overline{\cM}G_{0,k+1}({\cF^{(n)}}, \beta)_{\mu}$ be the $\textbf{S}$-fixed locus in $\overline{\cM}G_{0,k+1}({\cF^{(n)}}, \beta)$ consisting of objects for which the image of $\mathbb{P}^1$ lies on $Y^\mu \subset ({\cF^{(n)}})^{\bf{S}}$. 
Let $(F_{k_1,\beta_1}^{k_2,\beta_2})_\mu$ be the substack of $F_{k_1,\beta_1}^{k_2,\beta_2}$ taking the marked point to $Y^\mu$.

With these spaces, a short strategy of the proof is as follows.
Following the idea in \cite[Section \S 7.6]{CK0}, we express the function \eqref{eq:poly} for $G=J^{\bf S}_{\cF^{(n)}}$ as a sum of the integrations over $\overline{\cM}G_{0,k+1}({\cF^{(n)}}, \beta)_{\mu}$'s, see \eqref{PolyCond} below.
This seems reasonable because the latter (the RHS of \eqref{PolyCond}) is a sum of the integrations over $(F_{k_1,\beta_1}^{k_2,\beta_2})_\mu$'s by virtual localisation.
On the other hand it is equal to the former (the LHS of \eqref{PolyCond}) due to the isomorphism \eqref{hahaiso} and the expression \eqref{graphJform} of $J^{\bf S}_{\cF^{(n)}}$.
Then since $\overline{\cM}G_{0,k+1}({\cF^{(n)}}, \beta)_{\mu}$ is an $\bf S$-fixed locus, not a $\CC^*$-fixed locus, the latter is a polynomial in the $\CC^*$-equivariant parameter. In other words, it does not have a pole in $z=0$.

\smallskip

For each class $\beta$ and $i=1, ..., n$, we'd like to construct a $\mathbb{C}^* \times \bf{S}$-equivariant line bundle $E_{i,\beta}$ on $\overline{\cM}G_{0,k+1}(\cF^{(n)}, \beta)$ such that
\begin{equation} \label{Key:Equiv1}
E_{i,\beta}|_{F_{k_1 ,\beta_1}^{k_2, \beta_2}}\ =\ \text{ev}^*_\bullet (\det \cF_i^\vee \ot \pi^*\mathcal{O}_{Y}(1))\otimes \mathbb{C}_{\beta_2 (\det \cF_i^\vee \ot \pi^* \mathcal{O}_{Y}(1))},
\end{equation}
where $\mathbb{C}_{\star}$ denotes the $\CC^*$-representation of weight $\star$.
Using the $\bf S$-bundle $\det \cF_i^\vee \ot \pi^* \mathcal{O}_{Y}(1)$, we obtain an $\bf S$-morphism $\iota_i: \cF^{(n)} \ra \PP^{N_i}$ for some $N_i$ (which may not be an embedding).
It induces an $\bf S$-morphism 
\begin{align} \label{mora}
\overline{\cM}G_{0,k+1}(\cF^{(n)}, \beta)\ \longrightarrow\ \overline{\cM}G_{0,k+1}(\PP^{N_i}, \iota_{i*}\beta).
\end{align}
It automatically becomes a $\CC^*$-morphism as well by definition of the graph spaces \eqref{Gsp}.
By forgetting marked points, contracting all components except for the distinguished $\PP^1$ on the domain curve and replacing contracted points on $\PP^1$ with base points with the same degrees of the contracted components \cite[eq(3.2.3)]{CK0}, we obtain a morphism
\begin{align} \label{morb}
\overline{\cM}G_{0,k+1}(\PP^{N_i}, \iota_{i*}\beta)\ \longrightarrow \ \PP((\Sym^{\iota_{i*}\beta}H^0(\PP^1, \cO_{\PP^1}(1)))^{\oplus {N_i+1}}).
\end{align}
In \cite[eq(3.2.3)]{CK0}, the above morphism \eqref{morb} is written as 
$$
G_{0,k+1,\iota_{i*}\beta}(\PP^{N_i}) \ \longrightarrow \ QG_{0,k+1,\iota_{i*}\beta}(\PP^{N_i}).
$$ 
It is just a notational difference. They are exactly same \cite[Section \S 3.3, $2^{\mathrm{nd}}$ eq]{CK0}. The quasimap moduli space (in the sense of \cite{CKM})
$$
QG_{0,k+1,\iota_{i*}\beta}(\PP^{N_i}) \ =\ \PP((\Sym^{\iota_{i*}\beta}H^0(\PP^1, \cO_{\PP^1}(1)))^{\oplus {N_i+1}})
$$
is a compactification of a space of degree $\iota_{i*}\beta \in \ZZ$ morphisms from $\PP^1$ to $\PP^{N_i}$ by allowing base points.
The line bundle $E_{i,\beta}$ is defined as the pullback of the dual tautological bundle of the projective space by the map $\eqref{morb} \circ \eqref{mora}$.
See \cite[Section \S 3.3]{CK0} for more detailed construction. Equation \eqref{Key:Equiv1} is explained in \cite[eq(5.2.1)]{CK0}.

Using $\pi^* \left( \left(\det \cF_i^\vee\right)|_{Y^\mu} \ot \mathcal{O}_{Y}(1) \right)$, we can construct a $\mathbb{C}^* \times \bf{S}$-line bundle $E_{i,\beta,\mu}$ on $\overline{\cM}G_{0,k+1}(\cF^{(n)}, \beta)$ in the same way. Indeed, the construction comes through the forgetful (and stabilisation) morphism
$$
\overline{\cM}G_{0,k+1}(\cF^{(n)}, \beta)\ \ra\ \overline{\cM}G_{0,k+1}(Y, \pi_*\beta).
$$
Then the restriction of $E_{i,\beta,\mu}$ is
\begin{equation} \label{Key:Equiv2}
E_{i,\beta,\mu}|_{F_{k_1 ,\beta_1}^{k_2, \beta_2}}\ =\ \text{ev}^*_\bullet \pi^*(\mu^* \det \cF_i^\vee \ot \mathcal{O}_{Y}(1))\otimes \mathbb{C}_{\pi_*\beta_2 (\mu^* \det \cF_i^\vee \ot \mathcal{O}_{Y}(1))}.
\end{equation}

Then letting $\cE_{i,\beta, \mu} := E_{i, \beta} \ot E^\vee_{i,\beta,\mu}$, we obtain by \eqref{Key:Equiv1}, \eqref{Key:Equiv2} that
\begin{equation} \label{Key:Equiv}
\cE_{i,\beta,\mu}|_{F_{k_1 ,\beta_1}^{k_2, \beta_2} \cap \overline{\cM}G_{0,k+1}({\cF^{(n)}}, \beta)_{\mu}} = \mathbb{C}_{\beta_2 (\det \cF_i^\vee) +\pi_*\beta_2( \mu^*\det \cF_i )}  =  \CC_{\beta_2(E_{i,\mu})}.
\end{equation}
Let 
$$
N^{vir}_\mu\ :=\ N^{vir}_{\overline{\cM}G_{0,k+1}({\cF^{(n)}}, \beta)_{\mu}/\overline{\cM}G_{0,k+1}({\cF^{(n)}}, \beta)}
$$ 
be the virtual normal bundle with respect to the $\bf S$-action. 
By abuse of notation, we set $q^\beta := q^dQ^D$ for $\beta=(d,D)$, ${\bf t} := {\bf t} + \pi^*u$ and consider
{
\medmuskip=-1mu
\thinmuskip=-1mu
\thickmuskip=-1mu
\begin{align*} 
 Z_\mu &:=  \sum_{k,\beta \geq 0} \frac{q^\beta}{k!}  \frac{[\overline{\cM}G_{0,k+1}({\cF^{(n)}},\  \beta)_{\mu}]^{vir}}{e_{\CC^* \times \textbf{S}}(N^{vir}_\mu)} \ \cap\ e^{\sum_i c_1(\cE_{i,\beta, \mu} )y_i}  \prod_{a=1}^k \text{ev}_{a}^* (\textbf{t}) \text{ev}_{k+1}^*(p_0\pi^* \delta_j),
\end{align*}}
\par \noindent
where $p_0 := e_{\CC^*}(\cO(1) \ot \CC_1) \in H^*_{\CC^*}(\PP^1)$, so $p_0 |_0 = z$, and $p_0 |_\infty = 0.$
Recall from Section \S \ref{subsub:poly} that $\{\delta_i\}$ is a basis for $H^*(Y, \QQ)$ corresponding to the variables $\{u_i\}$.
We observe that $Z_\mu$ has no poles in $z=0$ because each factor of the denominator has a non-trivial $\lambda$-term.

For the morphism $\text{pt}: \overline{\cM}G_{0,k}({\cF^{(n)}}, \beta)_{\mu} \ra \Spec \; \CC$ to the point, we prove that 
\begin{align}
\label{PolyCond}
(z\partial_j J^{{\bf{S}},\mu}_{\cF^{(n)}}(z,q), J^{{\bf{S}},\mu}_{\cF^{(n)}}(-z,qe^{-z \sum_i y_iE_i }) )_Y\ =\  (\text{pt})_*Z_\mu,
\end{align}
using virtual $\CC^*$-localisation \cite{GP}. 
The fixed loci is a disjoint union of $(F_{k_1, \beta_1}^{k_2,\beta_2})_\mu$'s.
By abuse of notation, we denote by $N^{vir}$ the virtual normal bundle to $(F_{k_1, \beta_1}^{k_2,\beta_2})_\mu$ in $\overline{\cM}G_{0,k+1}({\cF^{(n)}}, \beta)$. Then we obtain 
{
\medmuskip=-1mu
\thinmuskip=-1mu
\thickmuskip=-1mu
\begin{align*} 
 Z_\mu &=   \sum_{k,\beta \geq 0}  \sum_{\tiny{\begin{array}{c}\tiny{k_1 +k_2 =k}\\ \tiny{\beta_1 +\beta_2=\beta}\end{array}}}  \frac{q^\beta}{k!}   \frac{[(F_{k_1 ,\beta_1}^{k_2, \beta_2})_\mu]^{vir}}{e_{\CC^* \times \textbf{S}}(N^{vir})} \ \cap\ e^{\sum_i c_1(\cE_{i,\beta,\mu} )y_i} \prod_{a=1}^k \text{ev}_{a}^* (\textbf{t}) \text{ev}_{k+1}^*(p_0\pi^* \delta_j).
\end{align*}}
\par \noindent
By \eqref{Key:Equiv} and the projection formula, we obtain $\mathrm{pt}_*$ of each term
\begin{align*}
 & (\text{pt})_*  \frac{q^\beta}{k!}   \frac{[(F_{k_1 ,\beta_1}^{k_2, \beta_2})_\mu]^{vir}}{e_{\CC^* \times \textbf{S} }(N^{vir})}  \cap e^{\sum_i c_1(\cE_{i,\beta, \mu} )y_i} \prod_{a=1}^k \text{ev}_{a}^* (\textbf{t})  \text{ev}_{k+1}^*(p_0\pi^* \delta_j)  \\
 &=  \frac{q^\beta}{k!} e^{ z\sum_i y_i \beta_2(E_{i,\mu}) } \int_{[F_{k_1 ,\beta_1}^{k_2, \beta_2}]^{vir}} \frac{(\text{ev}_\bullet)^*( PD(Y_\mu))}{e_{\mathbb{C}^* \times {\bf S}}(N^{k_2 ,\beta_2}_{k_1, \beta_1})} \prod_{a=1}^k \text{ev}_{a}^* (\textbf{t})  \text{ev}_{k+1}^*(p_0\pi^* \delta_j)  , \\
\end{align*}
where $N^{k_2,\beta_2}_{k_1,\beta_1} := N^{vir}_{F^{k_2,\beta_2}_{k_1,\beta_1}/\overline{\cM}G_{0,k+1}({\cF^{(n)}}, \beta)}$. 
$PD$ stands for `Poincar\'e dual'. 
Letting $N^{vir}_{k,\beta} := N^{vir}_{F_{k,\beta}/\overline{\cM}G_{0,k}({\cF^{(n)}}, \beta)}$ and $N_{k,\beta}^{'vir}$ be the bundle $N^{vir}_{k,\beta}$ with the inverse $\CC^*$-action, the integration on the right hand side is
\begin{align*} 
&\int_{[F_{k_1 ,\beta_1}^{k_2, \beta_2}]^{vir}} \frac{(\text{ev}_\bullet)^*( PD(Y_\mu))  }{e_{\mathbb{C}^* \times {\bf S}}(N_{k_1, \beta_1}^{k_2 ,\beta_2})} \prod_{a=1}^{k} \text{ev}^*_{a} (\textbf{t})  \text{ev}_{k+1}^*(p_0\pi^* \delta_j)  \\
& = \binom{k}{k_1} \int_{[F_{k_1 ,\beta_1}]^{vir}} \frac{z (\text{ev}_{\bullet})^*(\pi^*(\delta^s ) PD(Y_\mu)) \prod_{a=1}^{k_1} \text{ev}^*_{a} (\textbf{t}) \text{ev}_{k+1}^*(\pi^* \delta_j) }{e_{\mathbb{C}^* \times {\bf S}}(N_{ k_1, \beta_1}^{vir})} \times \\
& \qquad \qquad \qquad \qquad  \int_{[F_{k_2 ,\beta_2}]^{vir}}\frac{(\text{ev}_{\bullet})^*(\pi^*(\delta_s )PD(Y_{\mu})) \prod_{a=1}^{k_2} \text{ev}^*_{a} (\textbf{t})}{e_{\mathbb{C}^* \times {\bf S}}(N_{ k_2, \beta_2}^{'vir})}. 
\end{align*}
This proves the identity \eqref{PolyCond} holds true. Hence the left hand side of \eqref{PolyCond} has no poles in $z=0$.
\end{proof}

Since we have seen $S^*_{\bf t}(z)(\gamma)$ lies on the Lagrangian cone which is spanned by the $J$-function, it seems reasonable that $S^*_{\bf t}(z)(\gamma)$ satisfies the characterisation properties as well.

\begin{Prop} \label{Prop:condS}
The function $S^*_{{\bf{t}}}(z)(\gamma)$ satisfies the recursion relation and the polynomiality condition. The initial condition holds if $\gamma|_{q=Q=0}=1$. 
\end{Prop}
\begin{proof}
It satisfies the recursion relation by \cite[Theorem 2]{B}.
For the polynomiality condition, we can do the same computation in the proof of Proposition \ref{Prop:polyJ}. For the initial condition, we use \eqref{aaa} 
$S^*_{\bf{t}}(z)(1) = J^{\bf{S}}_{\cF^{(n)}}$ and the reconstruction of $\cL^\mu_Y$ together with \cite[Theorem 2]{B}.
\end{proof}

\subsection{Proof of Theorem \ref{Char}} \label{ProofChar}

Using the same ideas as in \cite[Lemma 6.4.1]{CK0} we obtain the following lemma.

\begin{Lemma} \label{Uptoz2}
Let $G  \in  \cH \; [\![{\bf t}, u]\!]  \subset  H^*_{\bf S} (\cF^{(n)}, \QQ)(z) \ot_\QQ \QQ[\![q, Q, {\bf t}, u]\!]$ be a function satisfying $G|_{q=Q=0} =e^{({\bf t}+\pi^*u)/z}$.
Then there exists 
\begin{itemize}
\item $\tau({\bf t},u)\ \in\ H^*_{\bf S}(\cF^{(n)}, \QQ)[\![q,Q, {\bf t}, u]\!]$, and 
\item $P({\bf t},u, z)\ \in\ H^*_{\bf S}(\cF^{(n)}, \QQ)[z] \ot_\QQ \QQ[\![q,Q, {\bf t}, u]\!]$ 
\end{itemize}
such that
\begin{enumerate}
\item $\tau ({\bf t},u)\ =\ {\bf t} + \pi^*u +O(q,Q)$,
\item $P({\bf t},u, z)\ =\ 1 +O(q,Q)$,
\item 
$
G ({\bf t})- S^*_{\tau({\bf t},u)} (z) (P(\tau({\bf t},u),z))\ \in\ z^{-2}H^*_{\bf S}(\cF^{(n)}, \QQ)[\![z^{-1},q, Q, {\bf t}, u]\!]
$
using the $z^{-1}$-expansion at simple poles.
\end{enumerate}
\end{Lemma}

\begin{proof}
We first write 
$$P \ : =\ \sum_{(d,D)} q^dQ^D P_{(d, D)} \ \ \text{ and }\ \  \tau \ : =\ \sum_{(d,D)} q^dQ^D \tau_{(d, D)} $$
and find $P_{(d, D)}$ and $\tau_{(d, D)}$ inductively.
We will use an induction on the partial order 
$$
(d',D') < (d,D) \iff D'<D \ \ \text{ or }\ \ D'=D, \ d'<d.
$$
For $(d,D)=0$, we set $P_0: =1$, $\tau_0({\bf t}, u) :={\bf t} + \pi^*u$ so that they satisfy (1), (2).
(3) for the coefficient of $q^0Q^0$ follows from the asymptotic property
\begin{align} \label{asympt}
S^*_{{\bf{t}}+\pi^*u}(z)(\gamma)\ =\ e^{({\bf{t}} +\pi^*u)/z} \gamma|_{q=0} + O(q)
\end{align}
and the initial condition for $G$.

Now, we fix $(d,D)$. Suppose that $P_{(d',D')}$ and $\tau_{(d',D')}$ are determined for all $(d',D')<(d,D)$ such that
after $z^{-1}$-expansion, the equality
\begin{align} \label{InDuction}
G\ =\  S^*_{\tau} (P(\tau)) ~~\text{ mod }{z^{-2}}
\end{align}
holds in all coefficients of $q^{d'}Q^{D'}$ with $(d',D') < (d,D)$.
By comparing the coefficients of $q^dQ^D$ in both sides of \eqref{InDuction}, we can uniquely determine $P_{(d,D)}$ and $\tau_{(d,D)}$ as follows. 
Using the asymptotic property \eqref{asympt} we write the right-hand side of \eqref{InDuction}
$$S^*_{\tau}(P(\tau))\ =\ e^{\frac{\tau( {\bf t }, u)}{z}}P(\tau) +O(q,Q).$$
The coefficient of $q^dQ^D$ in $O(q,Q)$ modulo $z^{-2}$ is already determined by the induction assumption.
The coefficient of $q^dQ^D$ in $e^{\frac{\tau({\bf t}, u)}{z}}P(\tau)$ is
$$e^{\frac{\textbf{t} + \pi^*u }{z}} \left( P_{(d,D)}+ \frac{1}{z}\tau_{(d,D)} \right)  + (\text{(mod $z^{-2}$)-known terms by the induction}).$$
So, $P_{(d,D)}+ \frac{1}{z}\tau_{(d,D)}$ is uniquely determined up to mod $z^{-2}$.
It determines $P_{(d,D)}$ and $\tau_{(d,D)}$.
\end{proof}

By Proposition \ref{Prop:condS}, $S^*_\tau (P(\tau))$ satisfies the characterisation properties.
Now, we would like to prove Theorem \ref{Char}. We assume $G$ satisfies the properties as well.
Let $F^\mu := \mu^* S^*_\tau (P(\tau))$. 
Lemma \ref{Uptoz2} implies that $G^\mu = F^\mu$ modulo $z^{-2}$, after $z^{-1}$-expansion, for all $\mu$.
The next lemma tells us that they are actually equal \emph{as rational functions}.
It follows from the idea in \cite[Lemma 7.7.1]{CK0}.

\begin{Lemma}  \label{Uniqueness}
We obtain $G^\mu = F^\mu$ as rational functions for all $\mu$.
\end{Lemma}

\begin{proof}
We first write 
$$G^\mu := \sum_{({\bf l}=(l_i),{\bf m}=(m_j),d,D)} q^d Q^D \prod_i t_i^{l_i} \prod_j u_j^{m_j} \cdot G^\mu_{({\bf l}=(l_i),{\bf m}=(m_j),d,D)} .$$
Let $F^\mu_{({\bf l},{\bf m},d,D)}$ 
be a function similarly defined using $F^\mu$.
For the proof, we're going to use an induction on $({\bf l},{\bf m},d,D)$ for the colexicographic order,
\begin{align*}
 & ({\bf l}',{\bf m}',d',D') < ({\bf l},{\bf m},d,D) \\
 & \iff  D' < D, \text{ or } D'=D, \ d' < d, \text{ or } D'=D, \ d'=d, \ {\bf m}' < {\bf m}, \\
 & \ \ \ \ \ \ \ \  \ \text{or } D'=D, \ d'=d, \ {\bf m}' = {\bf m}, \ {\bf l}'<{\bf l}.
\end{align*}
We would like to recall ${\bf l} \in \ZZ_{\geq 0}^{n}$, ${\bf m} \in \ZZ_{\geq 0}^{|\rank H^*(Y)|}$, $d \in \ZZ^n$, and $D \in \text{Eff}(Y)$.
As the first step of the induction,
we have $G^\mu_{({\bf l},{\bf m},0,0)} = F^\mu_{({\bf l},{\bf m},0,0)}$ for all ${\bf l}$ and ${\bf m}$ since both $G^\mu$ and $F^\mu$ follow the initial condition.

Now we fix $({\bf l},{\bf m},d,D)$ with $(d,D) \neq (0,0)$. 
Suppose that 
$$G^\mu_{({\bf l}',{\bf m}',d',D')}\ =\ F^\mu_{({\bf l}',{\bf m}',d',D')}$$
for all $({\bf l}',{\bf m}',d',D') < ({\bf l},{\bf m},d,D)$ and for all $\mu$. 
We will show that
\begin{align}  \label{WanTtoShoW}
G_{({\bf l},{\bf m},d,D)}^{\mu}\ =\ F_{({\bf l},{\bf m},d,D)}^\mu 
\end{align}
in several steps.

\subsubsection*{Step1}
For any ${\bf l}'$ and ${\bf m}'$, we have
\begin{align} \label{NoPoles}
G_{({\bf l}',{\bf m}',d,D)}^{\mu}-F_{({\bf l}',{\bf m}',d,D)}^\mu \ \in\ H^*_{\bf S} (Y, \QQ)[z, z^{-1}] 
\end{align}
for all $\mu$ as rational functions in $z$.
Here is a proof.
Since
$$
({\bf l}',{\bf m}',d-kd_{\mu,\nu},D)\ <\ ({\bf l},{\bf m},d,D)
$$
for any $k>0$ and $\nu$, we have
\begin{align} \label{ReCurSion}
\text{Res}_{z=-\frac{\chi_{\mu,\nu}}{k}} G^\mu_{({\bf l}',{\bf m}',d,D)}(z) dkz & = \frac{q^{kd_{\mu,\nu}}}{e_{\bf{S}}(N^{vir}_{\mu,\nu,k})} \mu^* \mu_* G^\nu_{({\bf l}',{\bf m}',d-kd_{\mu,\nu},D)} (-\frac{\chi_{\mu, \nu}}{k}) \nonumber \\
& =\frac{q^{kd_{\mu,\nu}}}{e_{\bf{S}}(N^{vir}_{\mu,\nu,k})} \mu^* \mu_* F^\nu_{({\bf l}',{\bf m}',d-kd_{\mu,\nu},D)} (-\frac{\chi_{\mu, \nu}}{k})  \\
& = \text{Res}_{z=-\frac{\chi_{\mu,\nu}}{k}} F^\mu_{({\bf l}',{\bf m}',d,D)}(z) dkz, \nonumber
\end{align}
where the first and third equalities follow from the recursion relation for $G^\mu$ and $F^\mu$, respectively, and the second equality comes from the induction assumption.
Then \eqref{NoPoles} is obtained by the pole conditions in the recursion relation and \eqref{ReCurSion}.

\subsubsection*{Step2}
By \eqref{NoPoles} and Lemma \ref{Uptoz2}, we obtain
\begin{align} \label{NoPoles1}
G_{({\bf l}',{\bf m}',d,D)}^{\mu}-F_{({\bf l}',{\bf m}',d,D)}^\mu \ \in\ z^{-2}H^*_{\bf S} (Y, \QQ)[z^{-1}] 
\end{align}
for any ${\bf l}'$, ${\bf m}'$, and $\mu$.

\subsubsection*{Step3}

Recall that $\{ \delta_j \}$ is a basis for $H^*(Y, \QQ)$. 
Let $\{\delta^j\}$ be its dual basis. 
We write $G^\mu_{({\bf l}, {\bf m}, d, D)}$ as
\begin{align*}
G^\mu_{({\bf l}, {\bf m}, d, D)}\ =\ \sum_{ j} G^{\mu,j}_{({\bf l},{\bf m},d,D)} \delta^j
\end{align*}
and consider the similar expression for $F^\mu_{({\bf l}, {\bf m}, d, D)}$.
Let $\Delta^{\mu, j}$ be the difference
{
\medmuskip=-1mu
\thinmuskip=-1mu
\thickmuskip=-1mu
\begin{align*}
(z\partial_{u_j}  G^\mu(z,q) , \ G^\mu(-z, qe^{-z\sum_i y_i E_{i,\mu} } ))_Y \ -\ (z\partial_{u_j}  F^\mu(z,q) , \ F^\mu(-z, qe^{-z\sum_i y_i E_{i,\mu} } ))_Y.
\end{align*}}
\par \noindent 
Let $\Delta^{\mu, j}_{({\bf l},{\bf m},d,D)}$ be the coefficient of degree $({\bf l},{\bf m},d,D)$-term in $\Delta^{\mu, j}$. 
Let ${\bf e}_j:=(0,...,0,1,0,..., 0) \in \ZZ_{\geq 0}^{|\rank H^*(Y)|}$ be the $j^{\mathrm{th}}$ standard basis. 
Then by the induction assumption, we can check that 
\begin{align} \label{Delta:Term1}
\Delta^{\mu,j}_{({\bf l},{\bf m},d,D)}  = \ & z(m_j +1)  (G_{({\bf l},{\bf m}+ {\bf e}_j,d,D)}^{\mu,0}-F_{({\bf l},{\bf m}+{\bf e}_j,d,D)}^{\mu,0})  
\\ & -(m_j +1)\sum_{j' \neq j}  (G_{({\bf l},{\bf m}+ {\bf e}_j-{\bf e}_{j'},d,D)}^{\mu,j'}-F_{({\bf l},{\bf m}+{\bf e}_j- {\bf e}_{j'},d,D)}^{\mu, j'})  \nonumber
\\ & -m_j  (G_{({\bf l},{\bf m},d,D)}^{\mu, j}-F_{({\bf l},{\bf m},d,D)}^{\mu, j}) \nonumber
\\ & +e^{-z\sum_i y_i (d_i - (d_{D,\mu})_i )}  (G_{({\bf l},{\bf m},d,D)}^{\mu, j}(-z)-F_{({\bf l},{\bf m},d,D)}^{\mu, j}(-z))  .\nonumber
\end{align}
By \eqref{NoPoles1}, we may write
\begin{align} \label{SubstiTute1}
G_{({\bf l},{\bf m},d,D)}^{\mu, j}-F_{({\bf l},{\bf m},d,D)}^{\mu, j}  \ =\ z^{-2}\left(A^{j,{\bf l},{\bf m}}_0 +\frac{A^{j,{\bf l},{\bf m}}_1}{z} + \cdots \right),
\end{align}
where $A^{j,{\bf l},{\bf m}}_0,A^{j,{\bf l},{\bf m}}_1, ... \in 
\QQ(\lambda)$ are uniquely determined $\lambda$-functions. 
Putting \eqref{SubstiTute1} into \eqref{Delta:Term1}, we observe that the coefficient of $z^{-1}$ is
\begin{align} \label{z1:term}
\sum_{s>0} \frac{1}{s!} \; A^{j,{\bf l},{\bf m}}_{s-1} \left(\sum_i y_i \left( d_i - (d_{D,\mu})_i \right) \right)^{s} + A^{0,{\bf l}, {\bf m}+{\bf e}_j}_0 .
\end{align}
The polynomiality conditions for $G^\mu$ and $F^\mu$ says \eqref{z1:term} must be zero for every $j$.
Hence 
we have $A^{j,{\bf l},{\bf m}}_s=0$ for all $s$ when $d  \neq d_{D, \mu}$.

\subsubsection*{Step4}
When $d  = d_{D, \mu}$, \eqref{WanTtoShoW} follows from the initial condition for both $G$ and $F$.
Since both $G^\mu_Y(-z)$ and $F^\mu_Y(-z)$ lie on $-z^{-1}\cL^\mu_Y$, and they are equal to each other up to $z^{-2}H^*_{\bf S} (Y, \QQ)[z^{-1}] \ot_\QQ \QQ[\![Q, {\bf t}, u]\!]$, we have $G^\mu_Y(-z)=F^\mu_Y(-z)$ by the reconstruction of $\cL^\mu_Y$ \cite[Proposition 1, or Section 8]{CG}.
\end{proof}

\begin{proof}[Proof of Theorem \ref{Char}]
Lemma \ref{Uniqueness} tells us that $G=S^*_\tau (P(\tau))$ where the right-hand side lies on $-z^{-1}\cL ag^{\bf S}_{\cF^{(n)}}$. 
This proves Theorem \ref{Char}.
\end{proof}

\section{Recursion relation for $I$-function} \label{RecurSion}

In this section, we would like to prove the following proposition.
\begin{Prop} \label{Prop:recI1}
For each fixed locus $\mu$, $I^{{\bf{S}},\mu}_{\cF^{(n)}}:=\mu^*I^{\bf S}_{\cF^{(n)}} $ has simple poles at $z=-\chi_{\mu, \nu}/k$ for $k \in \NN$, with residues
\begin{align} \label{RECURSIONaa}
\text{Res}_{z=-\frac{\chi_{\mu,\nu}}{k}} I^{{\bf{S}},\mu}_{\cF^{(n)}}(z) dkz\ =\ \frac{q^{kd_{\mu,\nu}}}{e_{\bf{S}}(N^{vir}_{\mu,\nu,k})} e_{\bf{S}}(N_{Y^\mu / \cF^{(n)}}) I^{{\bf{S}},\nu}_{\cF^{(n)}} (-\frac{\chi_{\mu, \nu}}{k}).
\end{align}
\end{Prop}
Note that Proposition \ref{Prop:recI1} does not guarantee that these are all simple poles.

\smallskip

Before providing a proof, we introduce some terminologies in combinatorics.
Each fixed locus $\mu: Y \rightarrow {\cF^{(n)}}$ is assigned to $I^\mu = (I^\mu_0, I^\mu_1, ..., I^{\mu}_{n}, I^\mu_{n+1})$ an ascending chain of subsets of $[r] := \{1,2,...,r\}$, 
$$\emptyset \; =:\; I^\mu_0\; \subset\; I_1^\mu\; \subset\; I^\mu _2\; \subset\; \cdots\; \subset\; I^\mu_n\; \subset\; I^\mu_{n+1 }\; :=\; [r], \ \ |I^\mu _i| \; =\; r_i.$$
Two different fixed loci $\mu$ and $\nu$ are connected by one-dimensional orbit if there are integers $\alpha_{\mu,\nu}, \beta_{\mu,\nu} \in [r]$ such that $I^\mu$ is replaced with $I^\nu$ by exchanging $\alpha_{\mu,\nu}$ and $\beta_{\mu, \nu}$.
Let $n^0_{\mu,\nu}, n^1_{\mu,\nu}$, $0 \leq n^0_{\mu,\nu} \leq n^1_{\mu,\nu} \leq n$, be indices determined by
\begin{align*}
n^0_{\mu,\nu} &\ :=\ \max \{i  :  | I^\mu_i \cap \{\alpha_{\mu,\nu} , \beta_{\mu,\nu}\} |\; =\; 0 \}, \\
n^1_{\mu,\nu} &\ :=\ \max \{i  :  | I^\mu_i \cap \{\alpha_{\mu,\nu} , \beta_{\mu,\nu}\} |\; =\; 1 \} .
\end{align*}
Since $\mu \neq \nu$, $n^0_{\mu,\nu} < n^1_{\mu,\nu}$. 
We distinguish $\alpha_{\mu,\nu}$ and $\beta_{\mu,\nu}$ by the property
\begin{align} \label{kakbconv}
\alpha_{\mu,\nu}\ \in\ I^\mu_{n^1_{\mu,\nu}},\ \beta_{\mu,\nu}\ \notin\ I^\mu_{n^1_{\mu,\nu}}
\end{align}
(hence $\alpha_{\mu,\nu}= \beta_{\nu,\mu}$).

The family of one-dimensional orbits connecting $Y^\mu$ to $Y^\nu$ is isomorphic to $\PP(L_{\alpha_{\mu,\nu}} \oplus L_{\beta_{\mu,\nu}}) \subset \cF^{(n)}$.
Note that the convention \eqref{kakbconv} is equivalent to 
$$Y^\mu\ \cong\ \PP\left(L_{\alpha_{\mu,\nu}} \oplus 0 \right), \ \ Y^\nu\ \cong\ \PP\left( 0 \oplus L_{\beta_{\mu,\nu}} \right).$$
On the corresponding component $Y \hookrightarrow \overline{\cM}_{0,2}(\cF^{(n)}, kd_{\mu,\nu})^{\bf{S}}$, 
the universal curve is also isomorphic to $\PP(L_{\alpha_{\mu,\nu}} \oplus L_{\beta_{\mu,\nu}})$ and the universal morphism from the universal curve $\PP(L_{\alpha_{\mu,\nu}} \oplus L_{\beta_{\mu,\nu}})$ to $\PP(L_{\alpha_{\mu,\nu}} \oplus L_{\beta_{\mu,\nu}}) \subset \cF^{(n)}$ is the $k$-covering morphism.
We denote it by $f_k$.

\begin{proof}[Proof of Proposition \ref{Prop:recI1}]

We divide the proof in several steps.

\subsubsection*{Step 1}
At the fist step, we'd like to check that $I^{{\bf S},\mu}_{\cF^{(n)}}$ has simple poles at $z=-\chi_{\mu,\nu}/k$, $k=1,2,3,...$.

The Atiyah-Bott localisation theorem tells us that for an equivariant class $K \in H^2(\PP(L_{\alpha_{\mu,\nu}} \oplus L_{\beta_{\mu,\nu}}), \QQ)$, we obtain
\begin{align} \label{ProJ:for}
d_{\mu,\nu}(K) \cdot 1 \ =\ p_*(K) \ =\ \frac{\mu^*K - \nu^*K}{\chi_{\mu,\nu}}\ \in\ H^0(Y, \ZZ)\ \cong\ \ZZ,
\end{align}
where $p: \PP(L_{\alpha_{\mu,\nu}} \oplus L_{\beta_{\mu,\nu}}) \ra Y$ is the projection morphism.

Recall that $H_{n+1,j} = -c_1(\pi^*L_j)$. 
For simplicity, we denote $H_j:= -c_1(\pi^*L_j)$ during the proof.
Then we obtain 
\begin{align} \label{Useful:Id}
0 &\ =\ \mu^*H_{\alpha_{\mu,\nu}} - \nu^*H_{\alpha_{\mu,\nu}}\ =\ \mu^*H_{\beta_{\mu,\nu}} - \nu^*H_{\beta_{\mu,\nu}} \\ \nonumber
\chi_{\mu,\nu} & \ =\ \mu^*H_{\alpha_{\mu,\nu}} - \nu^*H_{\beta_{\mu,\nu}}\ =\ \nu^*H_{\alpha_{\mu,\nu}} - \mu^*H_{\beta_{\mu,\nu}}
\end{align}
by applying $K=H_j$ and $\cS^\vee$ the dual of the tautological bundle over $\PP(L_{\alpha_{\mu,\nu}} \oplus L_{\beta_{\mu,\nu}})$.
Then by putting the equations \eqref{Useful:Id} to $\mu^*\eqref{flag}$, we observe that $I^{{\bf S},\mu}_{\cF^{(n)}}$ has simple poles at $z=-\chi_{\mu,\nu}/k$, $k=1,2,3,...$.

\subsubsection*{Step 2}
At the second step, we'd like to compute $ \frac{ e_{\bf{S}}(N^{vir}_{\mu,\nu,k}) }{e_{\bf{S}}(N_{Y^\mu /{\cF^{(n)}}})}$ explicitly.
From the definition of the virtual normal bundle $N^{vir}_{\mu,\nu,k}$ \cite{GP}, we have
{\medmuskip=-3mu
\thinmuskip=-3mu
\thickmuskip=-3mu
\begin{align} \label{zerothCal}
e_{\bf{S}}(N^{vir}_{\mu,\nu,k}) \ = \ &  \frac{e_{\bf{S}}(R^0p_*(f^*_kT_{{\cF^{(n)}}}|_{\PP(L_{\alpha_{\mu,\nu}} \oplus L_{\beta_{\mu,\nu}})}) ^{\text{mov}} ) }{e_{\bf{S}}(R^1p_*(f^*_kT_{{\cF^{(n)}}}|_{\PP(L_{\alpha_{\mu,\nu}} \oplus L_{\beta_{\mu,\nu}})}) ^{\text{mov}}  )  }\ \cdot \ \frac{-(\mu^*H_{\beta_{\mu,\nu}}  -\mu^*H_{\alpha_{\mu,\nu}})^2}{e_{\bf{S}}(R^0p_*( T_{p}))}  .
\end{align}}
\par \noindent
The first fractional term on the right-hand side is obtained by deformation-obstruction spaces of the universal morphism and the second fractional term on the right-hand side is obtained by automorphism spaces of the universal curve with two branch points.
Here, $T_p \cong \ker(T_{\PP(L_{\alpha_{\mu,\nu}} \oplus L_{\beta_{\mu,\nu}})} \ra p^* T_Y)$ is the relative tangent bundle of $p$.

We can compute $ \frac{ e_{\bf{S}}(N^{vir}_{\mu,\nu,k}) }{e_{\bf{S}}(N_{Y^\mu /{\cF^{(n)}}})}$ using \eqref{zerothCal}.
The restriction of the tangent bundle $T_{\cF^{(n)}}$ to $\PP(L_{\alpha_{\mu,\nu}} \oplus L_{\beta_{\mu,\nu}})$ becomes
{\medmuskip=-3mu
\thinmuskip=-3mu
\thickmuskip=-3mu
\nulldelimiterspace=-2pt
\scriptspace=-1.5pt
\begin{align} \label{Tangent}
T_{\cF^{(n)}} \ |_{\PP(L_{\alpha_{\mu,\nu}}\ \oplus L_{\beta_{\mu,\nu}}\ )} \ \cong & \ T_Y \ \oplus \  \oplus_{i \leq n^0_{\mu,\nu}} \Hom( \oplus_{j \in I^\mu_{i-1}}  L_j\ ,  \ \oplus_{j \in I^\mu_i \backslash I^\mu_{i-1}}  L_j\ ) \nonumber \\ 
& \oplus \ \Hom(\oplus _{j \in I^\mu_{r_{n^0_{\mu,\nu}}}} L_j \ ,  \ \oplus _{j \in I^\mu_{r_{n^0_{\mu,\nu}+1}} \backslash I^\mu_{r_{n^0_{\mu,\nu}} } , \ j \neq \alpha_{\mu,\nu}} L_j \  \oplus \  \cS ) \\ \nonumber
& \oplus \ \oplus_{n^0_{\mu,\nu} +1 <  i \leq n^1_{\mu,\nu}} \Hom(\oplus _{j \in I^\mu_{i } , j \neq \alpha_{\mu,\nu}} L_j \ \oplus \ \cS,\ \oplus _{j \in I^\mu_{i+1} \backslash I^\mu_{i } } \ L_j \ )\\ \nonumber
& \oplus \  \Hom(\oplus _{j \in I^\mu_{r_{n^1_{\mu,\nu}}}, j \neq \alpha_{\mu,\nu}} L_j \  \oplus  \ \cS  , \ \oplus _{j \in I^\mu_{r_{n^1_{\mu,\nu}+1}} \backslash I^\mu_{r_{n^1_{\mu,\nu}} } ,\ j \neq \beta_{\mu,\nu}} L_j \ \oplus \ \cS^\vee ) \\ \nonumber
& \oplus \ \oplus_{i > n^1_{\mu,\nu}+1} \ \Hom( \oplus_{j \in I^\mu_{i-1}}  L_j \ ,\  \oplus_{j \in I^\mu_i \backslash I^\mu_{i-1}}  L_j\  ) .
\end{align}}
\par \noindent
On the right-hand side, all bundles over $Y$ are considered to be pullbacks to $\PP(L_{\alpha_{\mu,\nu}} \oplus L_{\beta_{\mu,\nu}})$ along the projection morphism $p$.
Since the moving part occurs where $\cS$ or $\cS^\vee$ exist, we obtain 
\begin{align} \label{firstCal}
& \frac{e_{\bf{S}}(R^0p_*(f^*_kT_{{\cF^{(n)}}}|_{\PP(L_{\alpha_{\mu,\nu}} \oplus L_{\beta_{\mu,\nu}})})  ^{\text{mov}}  ) }{e_{\bf{S}}(R^1p_*(f^*_kT_{{\cF^{(n)}}}|_{\PP(L_{\alpha_{\mu,\nu}} \oplus L_{\beta_{\mu,\nu}})})  ^{\text{mov} }) } \\ 
 =&  \prod_{  \  l  \in  I^\mu_{ n^0_{\mu,\nu}}   }   \left( \prod_{s=1}^{k -1}(\mu^* H_{l}-\mu^* H_{\alpha_{\mu,\nu}} + s\frac{\chi_{\mu, \nu}}{k}) \right)^{-1} \nonumber \\
&    \prod_{  \  l'  \in ( I^\mu_{ n^1_{\mu,\nu}+1} \backslash I^\mu_{n^0_{\mu,\nu} +1} ) ,  l' \neq  \beta_{\mu,\nu} } \left(   \prod_{s=0}^{k }(\mu^* H_{\alpha_{\mu,\nu}}-\mu^* H_{l'}-s\frac{\chi_{\mu, \nu}}{k}) \right) \nonumber \\
&    \prod_{  \  l  \in  I^\mu_{ n^1_{\mu,\nu}}  ,  l \neq  \alpha_{\mu,\nu} }  \left(  \prod_{s=0}^{k }(\mu^* H_{l}-\mu^* H_{\beta_{\mu,\nu}}-s\frac{\chi_{\mu, \nu}}{k}) \right) \nonumber \\
& \prod_{s=0}^{2k}(\mu^*H_{\alpha_{\mu,\nu}}  -\mu^*H_{\beta_{\mu,\nu}}  -s\frac{\chi_{\mu,\nu}}{k} )
\nonumber
\end{align}
from \eqref{Tangent}.
Also we obtain  
\begin{align} \label{secondCal}
 \frac{-(\mu^*H_{\beta_{\mu,\nu}}  -\mu^*H_{\alpha_{\mu,\nu}})^2}{e_{\bf{S}}(R^0p_*( T_{p}))} \ =\ (\mu^*H_{\alpha_{\mu,\nu}}  -\mu^*H_{\beta_{\mu,\nu}}  -\chi_{\mu,\nu} )^{-1}
\end{align}
by a direct computation. Combining \eqref{Useful:Id}, \eqref{zerothCal}, \eqref{firstCal} with \eqref{secondCal}, we have
{\medmuskip=-3mu
\thinmuskip=-3mu
\thickmuskip=-3mu
\begin{align} \label{firstRes}
 \frac{ e_{\bf{S}}(N^{vir}_{\mu,\nu,k}) }{e_{\bf{S}}(N_{Y^\mu /{\cF^{(n)}}})} \  &=  \displaystyle       \prod_{  \  l  \in  I^\mu_{ n^0_{\mu,\nu}}   }   \left( \prod_{s=0}^{k -1}(\mu^* H_{l} \ -\ \mu^* H_{\alpha_{\mu,\nu}} + s\frac{\chi_{\mu, \nu}}{k}) \right)^{-1} \nonumber \\
&    \prod_{  \  l'  \in ( I^\mu_{ n^1_{\mu,\nu}+1} \backslash I^\mu_{n^0_{\mu,\nu} +1} ) ,  l' \neq  \beta_{\mu,\nu} } \left(   \prod_{s=1}^{k }(\mu^* H_{\alpha_{\mu,\nu}}- \ \mu^* H_{l'}\ -\ s\frac{\chi_{\mu, \nu}}{k}) \right)  \\
&    \prod_{  \  l  \in  I^\mu_{ n^1_{\mu,\nu}}  ,  l \neq  \alpha_{\mu,\nu} }  \left(  \prod_{s=1}^{k }(\mu^* H_{l}\ -\ \mu^* H_{\beta_{\mu,\nu}}-\ s\frac{\chi_{\mu, \nu}}{k}) \right) \nonumber \\
& \prod_{s=1}^{k-1}(\mu^*H_{\alpha_{\mu,\nu}}  -\ \mu^*H_{\beta_{\mu,\nu}}  -\ s\frac{\chi_{\mu,\nu}}{k} )
\prod_{s=0}^{k-1}(\mu^*H_{\beta_{\mu,\nu}}  -\ \mu^*H_{\alpha_{\mu,\nu}}  +\ s\frac{\chi_{\mu,\nu}}{k} ) . \nonumber
\end{align}}
\par \noindent
Using \eqref{Useful:Id}, one can see that the equation \eqref{firstRes} is equivalent to
\begin{align} \label{secondRes}
 \frac{ e_{\bf{S}}(N^{vir}_{\mu,\nu,k}) }{e_{\bf{S}}(N_{Y^\mu /{\cF^{(n)}}})}  &=  \displaystyle       \prod_{  \  l  \in  I^\nu_{ n^0_{\mu,\nu}}   }   \left( \prod_{s=1}^{k }(\nu^* H_{l}-\nu^* H_{\beta_{\mu,\nu}} - s\frac{\chi_{\mu, \nu}}{k}) \right)^{-1} \nonumber \\
&   \prod_{  \  l'  \in ( I^\nu_{ n^1_{\mu,\nu}+1} \backslash I^\nu_{n^0_{\mu,\nu} +1} ) ,  l' \neq  \alpha_{\mu,\nu} } \left(   \prod_{s=0}^{k -1}(\nu^* H_{\beta_{\mu,\nu}}-\nu^* H_{l'}+s\frac{\chi_{\mu, \nu}}{k}) \right)  \\
&   \prod_{  \  l  \in  I^\nu_{ n^1_{\mu,\nu}}  ,  l \neq  \beta_{\mu,\nu} }  \left(  \prod_{s=0}^{k -1}(\nu^* H_{l}-\nu^* H_{\alpha_{\mu,\nu}} + s\frac{\chi_{\mu, \nu}}{k}) \right) \nonumber \\
& \mathop{\prod_{s=0}}_{s \neq k}^{2k-1 }(\nu^*H_{\beta_{\mu,\nu}}  -\nu^*H_{\alpha_{\mu,\nu}}  +s\frac{\chi_{\mu,\nu}}{k} ).  \nonumber
\end{align}

\subsubsection*{Step 3}
At the third step, we'd like to show \eqref{RECURSIONaa} holds true.
From the formula \eqref{flag}, we observe that each coefficient of $I^{{\bf S}, \mu}_{\cF^{(n)}}(z)$ is
\begin{align} \label{Restriction}
  \sum_{\sum_l d_i^l =d_i} \prod_{i=1}^n & \left( \prod_{1 \leq l \neq l' \leq r_i}\frac{\prod_{s=-\infty}^{d_i^l -d_i^{l'}}(\mu^* H_{i,l}- \mu^* H_{i,l'}+sz) }{\prod_{s=-\infty}^{0}(\mu^*H_{i,l}-\mu^* H_{i,l'}+sz)}  \right.   \\ 
 & \left. \displaystyle  \prod_{1 \leq l \leq r_i,~ 1 \leq l' \leq r_{i+1}}\frac{\prod_{s=-\infty}^{0}( \mu^* H_{i,l}-\mu^* H_{i+1,l'}+sz)}{\prod_{s=-\infty}^{d_i^l -d_{i+1}^{l'}}(\mu^* H_{i,l}-\mu^* H_{i+1,l'}+sz)} \right) \nonumber \\ 
 =  \sum_{\sum_l d_i^l =d_i} \prod_{i=1}^n & \left( \prod_{1 \leq l \neq l' \leq r_i}\frac{\prod_{s=-\infty}^{d_i^l -d_i^{l'}}(\mu^* H_{\sigma_\mu(l)}- \mu^* H_{\sigma_\mu(l')}+sz) }{\prod_{s=-\infty}^{0}(\mu^*H_{\sigma_\mu(l)}-\mu^* H_{\sigma_\mu(l')}+sz)}  \right.  \nonumber \\ 
 & \left. \displaystyle  \prod_{1 \leq l \leq r_i,~ 1 \leq l' \leq r_{i+1}}\frac{\prod_{s=-\infty}^{0}( \mu^* H_{\sigma_\mu(l)}- \mu^* H_{\sigma_\mu(l')}+sz)}{\prod_{s=-\infty}^{d_i^l -d_{i+1}^{l'}}(\mu^* H_{\sigma_\mu(l)}- \mu^* H_{\sigma_\mu(l')}+sz)} \right) ,\nonumber \\ \nonumber
\end{align}
where $\sigma_\mu: [r] \ra [r] \in \fS_r$ is a permutation of $r$ such that $[r_{i+1}] \backslash [r_i]$ maps to $I^\mu_{i+1} \backslash I^\mu_i$, $i=1, ..., n$.
Note that $\sigma_\mu$ can be taken \emph{any such} permutation because of the symmetry of $d_i^l$ for each $i$.
We compute a residue of \eqref{Restriction} at $z=-\chi_{\mu,\nu}/k$.
It is enough to consider the case when the denominator of the right-hand side of \eqref{Restriction} contains a factor $(\mu^* H_{\alpha_{\mu,\nu}} - \mu^* H_{\beta_{\mu,\nu}} + kz)$.
Then the residue is
\begin{align} \label{bbbb}
\left(\eqref{Restriction} \cdot (\mu^* H_{\alpha_{\mu,\nu}} - \mu^* H_{\beta_{\mu,\nu}} + kz) \right) \left|_{z=-\chi_{\mu,\nu} /k} \right.
\end{align}

To verify \eqref{RECURSIONaa}, it is enough to show that \eqref{bbbb} is equal to
{\medmuskip=-2mu
\thinmuskip=-2mu
\thickmuskip=-2mu
\nulldelimiterspace=-1pt
\scriptspace=0pt
\begin{align} \label{residueresidue}
\frac{ e_{\bf{S}}(N^{vir}_{\mu,\nu,k}) }{e_{\bf{S}}(N_{Y^\mu /{\cF^{(n)}}})} \ & \sum_{\sum_l d_i^l =d'_i} \ \prod_{i=1}^n  \left( \prod_{1 \leq l \neq l' \leq r_i} \ \frac{\prod_{s=-\infty}^{d_i^l -d_i^{l'}}(\nu^* H_{\sigma_\nu(l)}\ -\ \nu^* H_{\sigma_\nu(l')}\ -\ s\chi_{\mu,\nu}/k) }{\prod_{s=-\infty}^{0}(\nu^*H_{\sigma_\nu(l)}\ -\ \nu^* H_{\sigma_\nu(l')}\ -\ s\chi_{\mu,\nu}/k)}  \right.  \nonumber \\  
 & \left.   \prod_{1 \leq l \leq r_i,~ 1 \leq l' \leq r_{i+1}}\ \frac{\prod_{s=-\infty}^{0}( \nu^* H_{\sigma_\nu(l)}\ -\ \nu^* H_{\sigma_\nu(l')}\ -\ s\chi_{\mu,\nu}/k)}{\prod_{s=-\infty}^{d_i^l -d_{i+1}^{l'}}(\nu^* H_{\sigma_\nu(l)}\ -\ \nu^* H_{\sigma_\nu(l')}\ -\ s\chi_{\mu,\nu}/k)} \right)\  ,  
\end{align}}
\par \noindent
where 
$$d'_i \ =\ \left\{
\begin{array}{cl}
d_i - k & \text{if }\ n^0_{\mu,\nu}<  i \leq n^1_{\mu,\nu} , \\
d_i & \text{if } \ i \leq n^0_{\mu,\nu} \ \text{ or }\ n^1_{\mu,\nu} < i,
\end{array}
\right.
$$
and $\sigma_\nu =\sigma_\mu \circ (\alpha_{\mu,\nu}, \beta_{\mu,\nu}) \in \fS_r$ is a permutation of $[r]$, because 
\begin{align*}
& \sum_{\sum_l d_i^l =d'_i} \prod_{i=1}^n  \left( \prod_{1 \leq l \neq l' \leq r_i}\frac{\prod_{s=-\infty}^{d_i^l -d_i^{l'}}(\nu^* H_{\sigma_\nu(l)}- \nu^* H_{\sigma_\nu(l')}-s\chi_{\mu,\nu}/k) }{\prod_{s=-\infty}^{0}(\nu^*H_{\sigma_\nu(l)}-\nu^* H_{\sigma_\nu(l')}-s\chi_{\mu,\nu}/k)}  \right.  \nonumber \\  
 & \left. \prod_{1 \leq l \leq r_i,~ 1 \leq l' \leq r_{i+1}}\frac{\prod_{s=-\infty}^{0}( \nu^* H_{\sigma_\nu(l)}- \nu^* H_{\sigma_\nu(l')}-s\chi_{\mu,\nu}/k)}{\prod_{s=-\infty}^{d_i^l -d_{i+1}^{l'}}(\nu^* H_{\sigma_\nu(l)}- \nu^* H_{\sigma_\nu(l')}-s\chi_{\mu,\nu}/k)} \right) 
\end{align*}
is the coefficient of $I^{{\bf S}, \nu}_{\cF^{(n)}}(-\chi_{\mu,\nu}/k)$ and we have an identity 
\begin{align*} 
e^{-\frac{\mu^*{\bf{t}}}{ \chi_{\mu,\nu} / k}} q^d e^{\sum_i t_i d_i} \ =\ q^{kd_{\mu,\nu}} e^{-\frac{\nu^*{\bf{t}}}{ \chi_{\mu,\nu} / k}} q^{d'} e^{\sum_i t_i d'_i}.
\end{align*} 
induced by \eqref{ProJ:for} and 
\begin{align*} 
d_{\mu,\nu}(\det \cF_i^\vee) \ =\ \left\{ 
\begin{array}{cl}
1 & \text{if }\ n^0_{\mu,\nu}<  i \leq n^1_{\mu,\nu}, \\
0 & \text{if }\ i \leq n^0_{\mu,\nu} \ \text{ or }\ n^1_{\mu,\nu} < i.
\end{array}
\right.
\end{align*}
Here, $(\alpha_{\mu,\nu}, \beta_{\mu,\nu})$ denotes the permutation exchanging $\alpha_{\mu,\nu}$ and $\beta_{\mu,\nu}$.
For any integers $m$ and $l \neq \alpha_{\mu,\nu}, \beta_{\mu,\nu}$, we obtain
{
\medmuskip=-1mu
\thinmuskip=-1mu
\thickmuskip=-1mu
\nulldelimiterspace=-1pt
\scriptspace=0pt
\begin{align} \label{uuu}
& \frac{\prod_{s=-\infty}^{m}( \mu^* H_{\alpha_{\mu,\nu}}- \mu^* H_{l} \ -\ s\chi_{\mu,\nu} / k)}{\prod_{s=-\infty}^{0}(\mu^* H_{\alpha_{\mu,\nu}}- \mu^* H_{l} \ -\ s\chi_{\mu,\nu} / k)} \\ \nonumber
& =\ \frac{\prod_{s=-\infty}^{m}( \nu^* H_{\beta_{\mu,\nu}} + \chi_{\mu,\nu} \ - \nu^* H_{l} - s\chi_{\mu,\nu} / k)}{\prod_{s=-\infty}^{0}(\nu^* H_{\beta_{\mu,\nu}} +\chi_{\mu,\nu} \ - \nu^* H_{l} - s\chi_{\mu,\nu} / k)} \\ \nonumber  
& = \ \frac{\prod_{s=-\infty}^{m-k}( \nu^* H_{\beta_{\mu,\nu}}- \nu^* H_{l} \ -\ s\chi_{\mu,\nu} / k)}{\prod_{s=-\infty}^{0}(\nu^* H_{\beta_{\mu,\nu}}- \nu^* H_{l} \ -\ s\chi_{\mu,\nu} / k)}  
\ \prod_{s=0}^{k -1}(\nu^* H_{\beta_{\mu,\nu}}-\nu^* H_{l}\ +\ s\frac{\chi_{\mu, \nu}}{k}\ ),
\end{align}}
\par \noindent
where the first equality comes from \eqref{Useful:Id}.
Doing the same computations (for $l=\beta_{\mu,\nu}$ and for the inverse of the most LHS of \eqref{uuu}), putting all them into \eqref{Restriction} and using \eqref{secondRes}, we obtain $\eqref{bbbb}=\eqref{residueresidue}$.
\end{proof}

\section{Quasimaps to GIT fiber bundles} \label{QuasiMapModuli}
In \cite{CKM}, the quasimap moduli spaces $Q_{g,k}(X, \beta)$ for a GIT quotient $X = W /\!\!/ \bf{G}$ 
are constructed, which provide another compactification of $\cM_{g,k}(X, \beta)$ \eqref{modulii}.
In \cite{CKt, CK0}, Ciocan-Fontanine and Kim showed that a certain generating function of invariants of `graph' quasimap spaces lies on $\cL ag_X$ and it is equal to Givental's $I$-function when $X$ is a (complete intersection in) toric variety or partial flag variety \cite{BCK1, CKt}.
In this section, we will construct quasimap spaces for GIT fiber bundles.
In Section \S \ref{PolyI}, we will show that $I^{\bf{S}}_{\cF^{(n)}}$ can be written in terms of `graph' quasimap spaces. Using this, we will show that it satisfies the characterisation properties so that we can conclude $I^{\bf S}_{\cF^{(n)}}$ lies on the cone by Theorem \ref{Char}.

\subsection{GIT quotients and its presentations} \label{Sect:pr}

Let $\textbf{G}$ be a linearly reductive algebraic group acting on an affine variety $V$. Suppose that $V$ has at worst locally complete intersection singularities.
Let $Y$ be a smooth projective variety and 
$\pi: E \rightarrow Y$ 
be a fiber bundle with the fiber $V$.
Suppose there exists a fiberwise $\textbf{G}$-action satisfying the following property. 
For each $y \in Y$, there is an affine neighborhood $U \subset Y$, and a $\textbf{G}$-equivariant isomorphism 
$$
\varphi :\pi^{-1}(U) \ :=\ E \times_Y U \rightarrow V \times U
$$ 
such that $p_2 \circ \varphi = \pi|_{\pi^{-1}(U)}$, where $p_2:V \times U \rightarrow U$ is the projection morphism. 

For a fixed character $\theta \in \chi(\textbf{G}):= \Hom({\bf G}, \CC^*)$, we define a $GIT$ $quotient$ 
$$E /\!\!/_{\theta}\textbf{G}\ :=\ \text{Proj} \bigoplus^{\infty}_{n=0}(R^0 \pi_{*}(E \times \CC_{\theta})^{\otimes n})^{\textbf{G}},$$
where $\mathbb{C}_{\theta}$ is the one-dimensional representation of $\textbf{G}$ determined by $\theta$. 
Then $E /\!\!/_{\theta} \textbf{G}$ is a fiber bundle over $Y$ with the fiber $V /\!\!/_{\theta} \textbf{G}$. 

Suppose that there is a morphism of varieties $\psi: Y \rightarrow \prod_{j=1}^r \mathbb{P}^{n_j -1}$
for some $r \in \mathbb{Z}_{>0}$ and $n_1, \cdots, n_r \in \mathbb{Z}_{>0}$, a $\textbf{T} :=(\mathbb{C}^*)^r$-action on $V$ which commutes with the $\bf{G}$-action on $V$, and a positive integer $m \in \mathbb{Z}_{>0}$,

\begin{align} \label{dotteddata}
\left(\psi: Y \rightarrow \prod_{j=1}^r \mathbb{P}^{n_j -1}, \ ({\bf{T}} \times {\bf{G}})\mathrm{-action} \ \mathrm{on}\ V, \ m \in \mathbb{Z}_{> 0}\right),
\end{align}

\noindent which satisfies the following five conditions,

\begin{enumerate}
\item $E$ is the pullback of a vector bundle {\footnotesize $[(\prod_{j=1}^r \mathbb{C}^{n_j}) \times V / \textbf{T}]$ on
$[(\prod_{j=1}^r \mathbb{C}^{n_j}) / \textbf{T}]$} under the composition morphism between stacks
$$
\ \ \ Y \  \xrightarrow{\ \psi\ } \ \prod_{j=1}^r \mathbb{P}^{n_j -1}\ \hookrightarrow\ \prod_{j=1}^r[\mathbb{C}^{n_j} / \mathbb{C}^*] \; \cong\; \left[ \left(\prod_{j=1}^r \mathbb{C}^{n_j}\right) / \textbf{T}\right].
$$
In other words, we have a fiber product
\begin{equation} 
\label{diaTar}
\xymatrix@R=5mm{
E \ar[r]\ar[d] & [(\prod_{j=1}^r \mathbb{C}^{n_j}) \times V / \textbf{T}] \ar[d] \\
Y \ar[r] & [(\prod_{j=1}^r \mathbb{C}^{n_j}) / \textbf{T}].
}
\end{equation}
To simplify notation, we let $W:= \prod_{j=1}^r \mathbb{C}^{n_j}$ and $\widetilde{V} := W \times V$.
\item The morphism $E \rightarrow [\widetilde{V} / \textbf{T}]$
in (\ref{diaTar}) is $\textbf{G}$-equivariant. 
\item $\Gamma(\cO_{\widetilde{V}})^{\textbf{T} \times \textbf{G}} \cong \mathbb{C}$ as $\mathbb{C}$-algebras.
\item $V^s(\textbf{G}, \theta)=V^{ss}(\textbf{G}, \theta) \neq \emptyset$ and it is non-singular. 
Moreover, the $\textbf{G}$-action on $V^s(\textbf{G}, \theta)$ is free.
\end{enumerate}

\noindent Before stating the fifth condition, recall that
$\prod_{j=1}^r \mathbb{P}^{n_j -1} \cong W /\!\!/_{\bf 1} \textbf{T}$,
where ${\bf 1} = (1, ..., 1) \in \chi(\textbf{T}) \cong \ZZ^r$. 
Let $\widetilde{\theta}:=m{\bf 1} + \theta \in \chi(\textbf{T})\oplus \chi(\textbf{G}) \cong \chi(\textbf{T}\times \textbf{G})$. The fifth condition is then

\begin{enumerate}
\setcounter{enumi}{4}
\item $$\widetilde{V}^{ss}(\textbf{T} \times \textbf{G}, \widetilde{\theta})\ =\ W^{ss}(\textbf{T},{\bf 1}) \times V^{ss}(\textbf{G}, \theta).$$
\end{enumerate}

\smallskip

The condition (4) guarantees that $E /\!\!/_\theta \textbf{G}$ is non-singular and an open substack of $[E/\textbf{G}]$. 
Simply, we write $E /\!\!/ \textbf{G}$ instead of $E /\!\!/_{\theta} \textbf{G}$. 
We will label both $E /\!\!/ \textbf{G} \rightarrow Y$ and $[E/\textbf{G}] \rightarrow Y$ by $\pi$ when the context is clear,  
$$\xymatrix@R=7mm{
E/\!\!/ \textbf{G} \ar@{^{(}->}[r]\ar[dr]_-\pi & [E/ \textbf{G}] \ar[d]^-\pi \\
& Y.
}$$
From conditions (1) and (2), we have a fiber diagram
\begin{align}
\label{maindia}
\xymatrix@R=7mm{
[E/ \textbf{G}] \ar[r]\ar[d]^-\pi & [\widetilde{V} / \textbf{T} \times \textbf{G}] \ar[r]\ar[d] & [V / \textbf{T} \times \textbf{G}] \ar[d]\\
Y \ar[r] & [W / \textbf{T}] \ar[r] & [\text{Spec}\mathbb{C}/\textbf{T}].
} 
\end{align}
Conditions (3) and (5) guarantee that the GIT quotient 
$$
\widetilde{V} /\!\!/_{\widetilde{\theta}} (\textbf{T} \times \textbf{G}) \ \cong \ [\widetilde{V}^{ss}(\textbf{T} \times \textbf{G}, \widetilde{\theta}) / \textbf{T} \times \textbf{G}]
$$
is a non-singular, projective, open substack of $[\widetilde{V}/\textbf{T} \times \textbf{G}]$.
Moreover,
$$
\xymatrix@R=6mm{
E /\!\!/ \textbf{G} \ar[r]\ar[d]^-\pi & \widetilde{V} /\!\!/_{\widetilde{\theta}} (\textbf{T} \times \textbf{G}) \ar[d] \\
Y \ar[r] &  W /\!\!/_{{\bf 1}} \textbf{T}
}$$
is a fiber diagram, hence $E /\!\!/ \textbf{G}$ is projective as well.

\begin{Def} \label{def:ref} A \emph{presentation} of $(E, {\bf{G}}, \theta)$ is a data \eqref{dotteddata} satisfying the conditions (1)-(5).
\end{Def}

\subsection{Presentation of ${\cF^{(n)}}$} \label{rep:F}

Recall that $L_j \rightarrow Y$, $j=1, \cdots, r$, are line bundles on $Y$, and $F:=\bigoplus_{j=1}^r L_j \rightarrow Y$ is the sum of $L_j$. 
We define
$$
E\ :=\ \bigoplus_{i=1}^{n-1} \text{Hom}(\mathcal{O}_Y^{\oplus r_i}, \mathcal{O}_Y^{\oplus r_{i+1}}) \oplus \text{Hom}(\mathcal{O}_Y^{\oplus r_n}, F)\ \rightarrow\ Y
$$ 
for given numbers $0=r_0<r_1<\cdots <r_n<r_{n+1}=r$
where $\text{Hom}(A,B)$ denotes the space $\Spec_{\cO_Y} \Sym \cH om_{\cO_Y} (A,B)^\vee$. 
We define
$\textbf{G}:=\prod_{i=1}^n GL_{r_i}(\mathbb{C})$-action on $E$ by
$$
E \times \textbf{G} \ \rightarrow\ E, \ \ \  ((B_i)_{i=1}^n \in E,\ (A_i)_{i=1}^n \in {\bf G})\ \mapsto\ (A_{i+1}^{-1} \cdot B_i \cdot A_i)_{i=1}^n
$$
where $A_{n+1} = \text{Id}_r.$ 
Let $\theta:=\text{det}^n \in \chi(\textbf{G})$. Then $V^s(\textbf{G}, \theta)=V^{ss}(\textbf{G}, \theta) \neq \emptyset$, and it is non-singular. 
Moreover, the $\textbf{G}$-action on $V^s(\textbf{G}, \theta)$ is free. 
The quotient
$E /\!\!/ \textbf{G}$ is then isomorphic to ${\cF^{(n)}}$.

Choose any ample line bundle $\mathcal{O}(1)$ on $Y$ and let
$$
E' \ :=\  \bigoplus_{i=1}^{n-1} \mathrm{Hom}(\mathcal{O}_Y^{\oplus r_i}, \mathcal{O}_Y^{\oplus r_{i+1}}) \oplus \mathrm{Hom}(\mathcal{O}_Y^{\oplus r_n}, F(l))
$$
be a twisting of $E$ for an integer $l>0$. Then a $\bf{G}$-action on $E'$ is defined by the same way as it acts on $E$.
\begin{Prop}
\label{flagrepr}
The quotients are isomorphic
$$
E' /\!\!/ {\bf{G}}\ \cong\ E /\!\!/ \bf{G}
$$
as fiber bundles. Moreover
$(E', {\bf{G}}, \theta)$ has a presentation (Definition \ref{def:ref}) if $l$ is large enough.
\end{Prop}

\begin{proof}
Take a positive integer $l$ such that $L_j' := L_j(l)$, $j=1, \cdots, r$, are generated by their global sections.
Then
$E /\!\!/_{\theta} \textbf{G} \cong E' /\!\!/_{\theta} \textbf{G}.$
Letting $n_j:=\text{dim}H^0(Y,L_j')$, we have a morphism
$$
\psi_j\; :\; Y\ \rightarrow\ \mathbb{P}(H^0(Y,L_j')^\vee)\; \cong\; \mathbb{P}^{n_j -1}$$
such that $\psi_j^*(\mathcal{O}_{\mathbb{P}^{n_j -1}}(1)) \cong L_j'$, and therefore, we obtain
$\psi: Y \rightarrow \prod_{j=1}^r \mathbb{P}^{n_j -1}.$
Define a $\textbf{T} := (\mathbb{C}^*)^r$-action on $V= \bigoplus_{i=1}^{n} \text{Hom}(\mathbb{C}^{r_i}, \mathbb{C}^{r_{i+1}})$ by 
$$
\textbf{T} \times V\; \rightarrow\; V,  \ ((t_1, \cdots , t_r), (b_i)_{i=1}^n)\; \mapsto\; (b_1, \cdots, b_{n-1}, \text{diag}(t_1, \cdots, t_r) \cdot b_n).
$$
The conditions (1) and (2) for $(E' , {\bf{G}}, \theta)$ are satisfied with the $\psi$ and $\textbf{T} \times \textbf{G}$-action on $V$ above.

Consider the following $\textbf{T} \times \mathbb{C}^*$-action on $V$, 
\begin{align} \label{V}
V = \bigoplus_{i=1}^{n-1} \text{Hom}(\mathbb{C}^{r_i}, \mathbb{C}^{r_{i+1}}) \oplus \text{Hom}(\mathbb{C}^{r_n}, \bigoplus_{j=1}^r \mathbb{C}_{e_j + \mathrm{id}_{\CC^*}}),
\end{align}
where $e_j$, $j=1,\cdots,r$, are standard basis for $\chi(\textbf{T}) \cong \mathbb{Z}^r$, and $\mathrm{id}_{\mathbb{C}^*} \in  \chi(\mathbb{C}^*) $.
Note that we consider the trivial action on $\CC$ without a subscript.
Let $u_j := e_j + \mathrm{id}_{\CC^*} \in \chi(\textbf{T} \times \mathbb{C}^*)$ for each $j$. 
Let
$$
{\bf 1} \; : = \; e_1 + \cdots + e_r \ \ \text{ and }\ \ \alpha\; =\; \sum_{i=1}^n r_i \cdot \mathrm{id}_{\CC^*}\ \in\ \chi(\textbf{T} \times \mathbb{C}^*).
$$ 
Note that letting 
$$\lambda\; :\; \mathbb{C}^* \ \hookrightarrow\ \textbf{G}, \ \ t \ \mapsto\ (t \cdot \text{Id}_{r_i})_{i=1}^n,$$
we have $\alpha = \theta \circ \lambda$.
Now, we can choose $m \in \mathbb{Z}_{>0}$ satisfying the following conditions.

\begin{itemize}
\item $m{\bf 1} + \alpha$ is in the interior of the cone generated by $e_1, \cdots, e_r, u_j$ for all $j=1, \cdots, r.$
\item $m{\bf 1} + \alpha$ is not in the cone generated by $e_1, \cdots, \hat{e_j}, \cdots, e_r, u_1, \cdots, u_r$ for all $j=1, \cdots, r.$
\end{itemize}

\noindent Then $\left(\psi, \ ({\bf{T}}\times {\bf G}) \ra \mathrm{Aut} (V), \ m \right)$ becomes a presentation of $(E', {\bf{G}}, \theta)$.
\end{proof}

\subsection{Stable quasimaps} \label{StableQmap}

Suppose that $(E,\textbf{G},\theta)$ has a presentation
$$
\left(\psi: Y\rightarrow \prod_{j=1}^r \mathbb{P}^{n_j -1} ,\ \textbf{T} \times \textbf{G} \rightarrow \text{Aut}(V) ,\ m \in \mathbb{Z}_{>0}\right).
$$
We choose a class
\begin{align} \label{dc}
\beta\; =\; (\beta',\beta_0) \ \in\ \text{Ker}\left(\; \text{Pic}(Y)^\vee \oplus \text{Pic}^{\textbf{T} \times \textbf{G}}(V)^\vee \; \rightarrow\; \chi(\textbf{T})^\vee \right),
\end{align}
where the morphism is defined by using (\ref{maindia}). Here $(-)^\vee = \text{Hom}_{\mathbb{Z}}(-,\ZZ)$.
Note that for ${\cF^{(n)}}$ (with the presentation in section \S \ref{rep:F}), choosing $\beta$ is equivalent to choosing a class in $H_2({\cF^{(n)}} , \ZZ)$.
For $g,k \in \mathbb{Z}_{\geq 0}$ with $2g+k \geq 2$, consider the following data,
\begin{align} \label{QuasiMap}
\left(  \phi:(C,p_1, \cdots, p_k)\; \rightarrow\; (C_0,p_1,\cdots,p_k),\ f:C\; \rightarrow\; Y,\ P, \ u  \right)
\end{align}
satisfying
\begin{itemize}
\item $(C,p_1, \cdots, p_k)$ and $(C_0,p_1, \cdots, p_k)$ are prestable $k$-pointed curves of genus $g$,
\item $\phi$ is the contraction of all rational tails on $C$. A \emph{rational tail} is a maximal (with respect to a partial order defined by inclusions) connected tree of rational curves with no marked points, attached to other components at only one node on $C$,
\item $f$ is of degree $\beta'$ and $f$ restricted to each irreducible component of a rational tail is non-constant, 
\item $P$ is a $\textbf{G}$-principal bundle on $C_0$, 
\item $u\; \in\; \Gamma\left(C_0, P_f \times_{(\textbf{T} \times \textbf{G})} V \right)$, 
\item $\beta_0 (L)= \text{deg}\left(u^*(P_f \times_{(\textbf{T} \times \textbf{G})} L)\right)$ for all $L \in \text{Pic}^{\textbf{T} \times \textbf{G}}(V)$.
\end{itemize}
Here, the $(\textbf{T} \times \textbf{G})$-principal bundle $P_f$ may be defined on $C_0$ by using $f$ as follows (See \cite{CK0}, \cite{MOP}).
For $1 \leq j_0 \leq r$, having 
$$
C\ \xrightarrow{\ f \ }\  Y\ \xrightarrow{\ \psi\ }\  \prod_{j=1}^r \mathbb{P}^{n_j -1} \  \longrightarrow\  \mathbb{P}^{n_{j_0} -1}
$$
is equivalent to the existence of a surjective morphism of sheaves on $C$
$$
\mathcal{O}_C^{\oplus n_{j_0}} \ \xrightarrow{\ \varphi_{j_0}\ } \ \mathcal{L}_{j_0} \ \longrightarrow\  0,
$$
where $\mathcal{L}_{j_0}$ is the pullback of $\mathcal{O}_{\mathbb{P}^{n_{j_0} -1}}(1)$ on $C$. 
Let $T_1 , \cdots, T_l$ be rational tails on $C$ and $\widetilde{C} $ be the closure of $C \smallsetminus (\cup_{i=1}^l T_i)$. Let $t_1 , \cdots , t_l$ be the corresponding points of $T_1, \cdots , T_l$ at $\widetilde{C}.$ Consider the following morphism of sheaves on $\widetilde{C}$
$$
\mathcal{O}_{\widetilde{C}}^{\oplus n_{j_0}} \xrightarrow{\ \varphi_{j_0} |_{\widetilde{C}}\ }  \mathcal{L}_{j_0} |_{\widetilde{C}} \ \hookrightarrow\ \displaystyle \mathcal{L}_{j_0} |_{\widetilde{C}} \otimes \mathcal{O}_{\widetilde{C}}\left(\sum_{i=1}^l \text{deg}(\mathcal{L}_{j_0} |_{T_i}) \cdot t_i\right).
$$
Since $\widetilde{C}$ is isomorphic to $C_0$ through the composition $\widetilde{C} \hookrightarrow C \rightarrow C_0,$ the above map gives rise to the morphism of stacks
$$
C_0\ \rightarrow\ [\mathbb{C}^{n_{j_0}} / \mathbb{C}^*]
$$
with degree $\deg \mathcal{L}_{j_0}$.
So, we have the morphism of stacks
$$
C_0\ \rightarrow\ \prod_{j=1}^r [\mathbb{C}^{n_j} / \mathbb{C}^*]\ \cong\ [W / \textbf{T}]
$$
with degree $d:=(\deg \mathcal{L}_1 , \cdots, \deg \mathcal{L}_r) \in \mathbb{Z}^r \cong   \text{Hom}(\text{Pic}^{\textbf{T}}W, \mathbb{Z}) .$
Note that $d$ is the image of $\beta'$ under
$$
\text{Pic}(Y)^\vee\ \rightarrow\ \text{Pic}\left(\prod_{j=1}^r \mathbb{P}^{n_j -1}\right)^\vee \cong  \ \text{Pic}^{\textbf{T}}(W)^\vee.$$
In conclusion, we have the $\textbf{T}$-principal bundle $Q_f$ on $C_0$, and the $\textbf{T}$-equivariant morphism
\begin{align}
Q_f \ \rightarrow\ W. \label{Qbun}
\end{align}
Now, define $P_f := Q_f \times_{C_0} P.$

We call the above data \eqref{QuasiMap} a {\em genus $g$, $k$-pointed quasimap with degree $\beta$}, or simply, a {\em quasimap with type $(g,k,\beta)$}. Now, we want to define a stability condition on quasimaps. 
For a quasimap with type $(g,k,\beta)$,
$$
\left(  \phi:(C,p_1, \cdots, p_k)\; \rightarrow\; (C_0,p_1,\cdots,p_k),\ f:C\; \rightarrow\; Y,\ P, \ u  \right),
$$
we have the morphism $P_f  \rightarrow W$ through \eqref{Qbun}, which defines the section 
$C_0 \rightarrow P_f \times_{(\textbf{T} \times \textbf{G})} W.$
Combining this with $u$, we have the section
$$
\widetilde{u}\; :\; C_0 \ \rightarrow \ P_f \times_{(\textbf{T} \times \textbf{G})} (W\times V).$$
Thus, we obtain a 
quasimap $(C_0, P_f, \widetilde{u})$ in the sense of \cite{CKM}.

\begin{Def}
A quasimap 
$$(\phi:(C,p_1, \cdots, p_k) \rightarrow (C_0,p_1,\cdots,p_k),~~f:C \rightarrow Y,~~ P, ~~u)$$
is $\theta$-\emph{prestable} if 
$
(C_0, P_f, \widetilde{u})
$
is $(\widetilde{\theta} = m{\bf 1} + \theta)$-prestable in the sense of \cite{CKM}.
It is $\theta$-\emph{stable} if for each irreducible component $C' \subset C_0$, $f|_{C'}$ is non-constant or $(C_0, P_f, \widetilde{u})|_{C'}$ is $0^+$-stable in the sense of \cite{CKM} with respect to $\widetilde{\theta}$.
\end{Def}

\begin{Def}
\label{isoDef}
An \emph{isomorphism} between two quasimaps 
$$
\left(  \phi:(C,p_1, \cdots, p_k)\; \rightarrow\; (C_0,p_1,\cdots,p_k),\ f:C\; \rightarrow\; Y,\ P, \ u  \right), \ \text{ and}
$$
$$
\left( \phi':(C',p'_1, \cdots, p'_k)\; \rightarrow\; (C'_0,p'_1,\cdots,p'_k),\ f':C'\; \rightarrow\; Y,\ P', \ u' \right)
$$
is a tuple of isomorphisms 
$$
\left( \; q:C\; \xrightarrow{ \sim  } \; C' ,\ \xi : P\; \xrightarrow{ \sim  } \; q_0^*P' \; \right)$$
such that 
$$
f' \circ q \; =\; f,\ q \circ p_a\; =\; p'_a, \  \text{ and }  \ q_0^* u' \;=\; (\zeta, \xi)\circ u,
$$
where
\begin{itemize}
\item $q_0: C_0 \xrightarrow{\sim} C'_0$ is an isomorphism between the contracted curves induced by $q$,
\item $\zeta$ is an isomorphism $Q_f \xrightarrow{\sim} q_0^*Q'_{f'}$ between the $\bf{T}$-principal bundles on $C_0$ induced by $f$ and $f'$ 
which commutes the diagram
$$
\xymatrix@R=6mm{
Q_f \ar[r] \ar[dr] & q_0^*Q'_{f'} \ar[d]\\
& W.
}
$$
\end{itemize}
Note that $\zeta$ is uniquely determined by the condition $f' \circ q  = f$.
\end{Def}

We denote the moduli space of (pre)stable quasimaps with type $(g,k,\beta)$ by 
$$
Q_{g,k}^{(pre)}(E /\!\!/ \textbf{G} ,\beta) .
$$
Then there is a morphism of moduli spaces of prestable quasimaps
$$Q^{pre}_{g,k}(E /\!\!/ \textbf{G},\beta) \ \rightarrow \ Q^{pre}_{g,k} \left(\widetilde{V} /\!\!/_{\widetilde{\theta}} (\textbf{T} \times \textbf{G}) , \beta''\right),$$
where $\beta'' := \deg (C_0,P_f,\widetilde{u})$. The latter space is in the sense of \cite{CKM}, 
which is an Artin stack, locally of finite type over $\mathbb{C}$. 
Since the stability condition is an open condition, $Q_{g,k}(E /\!\!/ \textbf{G},\beta)$ is an open substack of $Q^{pre}_{g,k}(E /\!\!/ \textbf{G},\beta)$.

\subsection{Quasimap moduli spaces}

Let
$\mathfrak{M}_{g,k}(Y,\beta')$
be the stack of genus $g$, $k$-pointed prestable maps to $Y$ with degree $\beta'$. 
It is an Artin stack, locally of finite type over $\mathbb{C}.$ 
Let $\mathfrak{M}'_{g,k}(Y,\beta')$ be the substack of $\mathfrak{M}_{g,k}(Y,\beta')$ whose objects are non-constant on each component of a rational tail.
$\mathfrak{M}'_{g,k}(Y,\beta')$ is an open substack of $\mathfrak{M}_{g,k}(Y,\beta')$; see \cite[Lemma 5.1]{O}.
There is a morphism of stacks
\begin{eqnarray*}
&Q^{pre}_{g,k}(E /\!\!/ \textbf{G}, \beta) \ \rightarrow\ \mathfrak{M}'_{g,k}(Y,\beta'), \\
&\left( \phi:(C,p_1, \cdots, p_k) \rightarrow (C_0,p_1,\cdots,p_k),\ f ,\ P,\ u \right) \ \mapsto \ \left((C,p_1, \cdots, p_k),\ f \right).
\end{eqnarray*}

We define a morphism of stacks 
\begin{eqnarray*}
&Q^{pre}_{g,k}\left(\widetilde{V} /\!\!/_{\widetilde{\theta}} (\textbf{T} \times \textbf{G}) , \beta'' \right) \ \rightarrow\ Q^{pre}_{g,k}(W /\!\!/_{\bf 1} \textbf{T}, d), \\
& \left((C_0,{\bf p}),\ P,\ u\right) \ \mapsto \ \left((C_0,{\bf p}),\ Q:= P/\textbf{G}, \ u': C_0 \rightarrow Q \times_{\textbf{T}} W\right),
\end{eqnarray*}
where $u'$ is the composition  
$$
C_0\ \xrightarrow{\ u\ }\ P \times_{\textbf{T} \times \textbf{G}} \widetilde{V} \ \longrightarrow\ P \times_{\textbf{T} \times \textbf{G}} W\ \cong\ Q \times_{\textbf{T}} W.
$$
Then we have 
$$
\xymatrix@R=6mm{
Q^{pre}_{g,k}(E/\!\!/ \textbf{G},\beta)\ar[r]\ar[d] & Q^{pre}_{g,k}(\widetilde{V} /\!\!/_{\widetilde{\theta}} (\textbf{T} \times \textbf{G}) , \beta'') \ar[d]\\
\mathfrak{M}'_{g,k}(Y,\beta') & Q^{pre}_{g,k}(W /\!\!/_{\bf 1} \textbf{T}, d).
}
$$
In order to define a morphism of stacks
$$
\mathfrak{M}'_{g,k}(Y,\beta') \ \rightarrow\ Q^{pre}_{g,k}(W /\!\!/_{\bf 1} \textbf{T}, d),
$$
which makes the above diagram commutative, we use
the morphism of contraction of rational tails
\begin{align} \label{lem2.4.1} 
\mathfrak{M}_{g,k} \ \rightarrow \ \mathfrak{M}_{g,k}.
\end{align}  
Here, $\mathfrak{M}_{g,k}$ denotes the moduli stack of genus $g$, $k$-pointed prestable curves.
For $((C,p_1,\cdots ,p_k), ~~f) \in \mathfrak{M}'_{g,k}(Y,\beta')$, we obtain the contraction of rational tails $(C_0,p_1,\cdots ,p_k) \in \mathfrak{M}_{g,k}$ using \eqref{lem2.4.1}. 
We now define a morphism
\begin{eqnarray*}
&\mathfrak{M}'_{g,k}(Y,\beta') \rightarrow Q^{pre}_{g,k}(W /\!\!/_{\bf 1} \textbf{T}, d), \\
&\left((C,p_1,\cdots ,p_k),\ f\right) \ \mapsto \ \left( (C_0,p_1,\cdots ,p_k),\ Q_f ,\ u' : C_0 \rightarrow Q_f \times_{\textbf{T}} W\right).
\end{eqnarray*}
Then the diagram
\begin{align}
\label{fiberProd}
\xymatrix@R=6mm{
Q^{pre}_{g,k}(E /\!\!/ \textbf{G},\beta)\ar[r]\ar[d] & Q^{pre}_{g,k}(\widetilde{V} /\!\!/_{\widetilde{\theta}} (\textbf{T} \times \textbf{G}) , \beta'') \ar[d]\\
\mathfrak{M}'_{g,k}(Y,\beta') \ar[r] & Q^{pre}_{g,k}(W /\!\!/_{\bf 1} \textbf{T}, d)
}
\end{align}
is commutative.

\begin{Lemma} \label{Forproper}
The above commutative diagram \eqref{fiberProd} is a fiber diagram.
Moreover, the space
\begin{align} \label{fQ}
Q^{pre}_{g,k}(E/\!\!/ {\bf{G}},\beta) \ \cong \ Q^{pre}_{g,k}(\widetilde{V} /\!\!/_{\widetilde{\theta}} ({\bf{T}} \times {\bf{G}}),\beta'') \times_{Q^{pre}_{g,k}(W /\!\!/_{\bf 1} {\bf{T}}, d)} \mathfrak{M}'_{g,k}(Y,\beta')
\end{align}
is an Artin stack, locally of finite type over $\CC$.
\end{Lemma}

\begin{proof}
There exists a canonical way to recover an object in $Q^{pre}_{g,k}(E/\!\!/ {\bf{G}},\beta)$ from a given one in $Q^{pre}_{g,k}(\widetilde{V} /\!\!/_{\widetilde{\theta}} ({\bf{T}} \times {\bf{G}}),\beta'') \times_{Q^{pre}_{g,k}(W /\!\!/_{\bf 1} {\bf{T}}, d)} \mathfrak{M}'_{g,k}(Y,\beta') $ which defines an inverse morphism
$$Q^{pre}_{g,k}(\widetilde{V} /\!\!/_{\widetilde{\theta}} ({\bf{T}} \times {\bf{G}}),\beta'') \times_{Q^{pre}_{g,k}(W /\!\!/_{\bf 1} {\bf{T}}, d)} \mathfrak{M}'_{g,k}(Y,\beta') \ \ra \ Q^{pre}_{g,k}(E/\!\!/ {\bf{G}},\beta).$$

If $\mathcal{X}$, $\mathcal{Y}$ (over $\mathcal{Z}$) and $\mathcal{Z}$ are Artin stacks, locally of finite type over $\mathbb{C}$, then so is $\mathcal{X} \times_{\mathcal{Z}} \mathcal{Y}$.
\end{proof}

\begin{Cor} 
$Q_{g,k}(E /\!\!/ {\bf{G}} ,\beta)$ is a DM stack, locally of finite type over $\mathbb{C}.$
\end{Cor}

\subsection{Properness of $Q_{g,k}(E /\!\!/ \textbf{G}, \beta)$}
For a degree class $\beta=(\beta', \beta_0) $ and a line bundle $\mathcal{M} = \mathcal{M}' \boxtimes \mathcal{M}_0$ on $E /\!\!/ \textbf{G}$ with $\mathcal{M}' \in \text{Pic}(Y)$, $\mathcal{M}_0 \in \text{Pic}^{\textbf{T} \times \textbf{G}}(V),$
we define 
$$
\beta(\mathcal{M})\ :=\ \beta'(\mathcal{M}') \ +\ \beta_0 (\mathcal{M}_0).
$$
Then it is well-defined -- if $\mathcal{M}' \boxtimes \mathcal{M}_0 \cong \cN' \boxtimes \cN_0$, then 
$$
\beta'(\mathcal{M}')\  -\ \beta'(\cN') \ = \ \beta_0 (\cN_0)\ -\ \beta_0 (\mathcal{M}_0) 
$$
by \eqref{dc}. Let $\mathcal{O}(1)$ be any ample line bundle on $Y$ and let 
\begin{align}
\label{lb}
\mathcal{L}\ :=\ \pi^*\left( \mathcal{O}(1) \otimes \psi^*\left(\boxtimes_{j=1}^r \mathcal{O}_{\mathbb{P}^{n_j -1}}(m)\right)\right)\; \otimes\; i^*[E \times \mathbb{C}_{\theta} / \textbf{G}],
\end{align}
where $i: E /\!\!/ \textbf{G} \hookrightarrow [E/\textbf{G}]$ is the open immersion.
Then $\cL$ is an ample line bundle on $E /\!\!/ \textbf{G}$.
By \cite[Lemma 3.2.1]{CKM}, a quasimap with type $(g,k,\beta)$
$$\left(\phi:(C,p_1, \cdots, p_k) \; \rightarrow \; (C_0,p_1,\cdots,p_k),\ f:C\; \rightarrow \; Y,\ P, \ u\right)$$
satisfies $\beta(\mathcal{L}) \geq 0$, and 
$$
\beta(\mathcal{L})=0\  \iff\ \beta =0\ \iff\ (C, f) \text{ and } (C_0, P_f, \widetilde{u}) \text{ are constant.}
$$
Therefore, the number of irreducible components of the underlying curve $C$ of stable quasimaps to $E /\!\!/ \textbf{G}$ is bounded, and so is the number of irreducible components of $C_0$.
The boundedness of $Q_{g,k}(E /\!\!/ \textbf{G}, \beta)$ then follows from \cite[Theorem 3.2.4]{CKM}.
Hence, it is of finite type over $\CC$.
Using the valuative criteria and the properness of moduli of stable maps or stable quasimaps \cite[Proposition 4.3.1]{CKM}, we obtain the following proposition.

\begin{Prop}
\label{proper} $Q_{g,k}(E /\!\!/ {\bf{G}} ,\beta)$ is proper over $\mathbb{C}$.
\end{Prop}

\begin{proof}
We will use the fiber product \eqref{fQ} of moduli of prestable (quasi)maps.
For the valuative criteria let $R$ be a DVR (discrete valuation ring) over $\mathbb{C}$ with the quotient field $K$. Let $\Delta :=\text{Spec}R$ and $0 \in \text{Spec}R$ be its unique closed point. Then $\Delta^0 := \Delta\backslash\{0\}=\text{Spec}K$. 

We first prove separatedness.
Let
$$\left( \left( (\mathcal{C}^i_0,p^i_1,\cdots,p^i_k \right),\ \mathcal{P}^i,\ u^i),\ \left((\mathcal{C}^i,p^i_1,\cdots,p^i_k),\ f^i\right) \right), \ i=1,2,$$
be two objects in 
$$
Q_{g,k}(E /\!\!/ \textbf{G}, \beta)(\Delta)\; \subset\; Q^{pre}_{g,k}(\tilde{V}/\!\!/_{\tilde{\theta}} (\textbf{T} \times \textbf{G}),\beta'') \times_{Q^{pre}_{g,k}(W /\!\!/_{\bf 1} \textbf{T}, d)} \mathfrak{M}'_{g,k}(Y,\beta')(\Delta)
$$ 
which are isomorphic over $\Delta^0.$ We prove that they are isomorphic over $\Delta$.
Let 
$$
\left((\mathcal{C}'^i_0,p'^i_1,\cdots,p'^i_k),\ \mathcal{P'}^i,\ u'^i \right)\ \in\ Q_{g,k}(\tilde{V}/\!\!/_{\tilde{\theta}} (\textbf{T} \times \textbf{G}),\beta'')
$$ 
be the stabilisation of $((\mathcal{C}^i_0,p^i_1,\cdots,p^i_k),~\mathcal{P}^i,~u^i)$.
Since $Q_{g,k}(\tilde{V}/\!\!/_{\tilde{\theta}} (\textbf{T} \times \textbf{G}),\beta'')$ is proper over $\mathbb{C}$ \cite[Proposition 4.3.1]{CKM} (and hence separated), they are isomorphic over $\Delta$.
So, we may assume that they are identical
$$
\left((\mathcal{C}'^1_0,p'^1_1,\cdots,p'^1_k),\ \mathcal{P'}^1,\ u'^1\right)\ =\ \left((\mathcal{C}'^2_0,p'^2_1,\cdots,p'^2_k),\ \mathcal{P'}^2,\ u'^2\right).
$$
Denote it by $((\mathcal{C}'_0,p'_1,\cdots,p'_k),\mathcal{P'},u').$
By adding additional marked points $r_1,\cdots,r_l$ on $\mathcal{C}'_0$ until 
$$
\left((\mathcal{C}^1,p^1_1,\cdots,p^1_k,r_1,\cdots,r_l),\ f^1\right) \ \text{ and }\  \left((\mathcal{C}^2,p^2_1,\cdots,p^2_k,r_1,\cdots,r_l),\ f^2\right)$$ 
are stable, they are isomorphic over $\Delta$ because $\overline{\cM}_{g,k+l}(Y,\beta')$ is proper.
So, 
$$\left((\mathcal{C}^1,p^1_1,\cdots,p^1_k),\ f^1\right)\ \text{ and }\ \left((\mathcal{C}^2,p^2_1,\cdots,p^2_k),\ f^2\right)$$ 
are isomorphic over $\Delta$ as well.
Now since $((\mathcal{C}^i_0,p^i_1,\cdots,p^i_k),~\mathcal{P}^i,~u^i)$ is the pullback of $((\mathcal{C}'_0,p'_1,\cdots,p'_k),\mathcal{P'},u')$, they are also isomorphic over $\Delta$.
Hence the separatedness is satisfied.

The completeness is a bit more complicated.
Let
\begin{align}
\label{data}
\left(\left((\mathcal{C}_0,p_1,\cdots,p_k),\ \mathcal{P},\ u\right),\ \left((\mathcal{C},p_1,\cdots,p_k),\ f \right)\right)
\end{align}
be a stable object in $Q^{pre}_{g,k}(\tilde{V}/\!\!/_{\tilde{\theta}} (\textbf{T} \times \textbf{G}),\beta'') \times_{Q^{pre}_{g,k}(W /\!\!/_{\bf 1} \textbf{T}, d)} \mathfrak{M}'_{g,k}(Y,\beta')(\Delta^0)$.
We would like to find its extension over $\Delta$ (by shrinking it if necessary).
Let $((\mathcal{C}'_0,p'_1,\cdots,p'_k),~\mathcal{P'},~u')$ be the stabilisation of $((\mathcal{C}_0,p_1,\cdots,p_k),\ \mathcal{P},\ u)$.
Since $Q_{g,k}(\tilde{V}/\!\!/_{\tilde{\theta}} (\textbf{T} \times \textbf{G}),\beta'')$ is proper over $\mathbb{C}$, we have its extension 
$$
\left((\bar{\mathcal{C}}'_0,\bar{p}'_1,\cdots,\bar{p}'_k),\ \bar{\mathcal{P}}',\ \bar{u}' \right)  \ \in \ Q_{g,k}(\tilde{V}/\!\!/_{\tilde{\theta}} (\textbf{T} \times \textbf{G}),\beta'')(\Delta).
$$
Hence, we have a family of prestable quasimaps 
$$
\left((\bar{\mathcal{C}}'_0,\bar{p}'_1,\cdots,\bar{p}'_k),\ \bar{\mathcal{Q}}'=\bar{\mathcal{P}}'/\textbf{G},\ \bar{u}'_{\textbf{T}}:\bar{\mathcal{C}}'_0 \rightarrow \bar{\mathcal{Q}}' \times_{\textbf{T}} W\right)  \ \in \ Q^{pre}_{g,k}(W/\!\!/_{\bf 1} \textbf{T} , d)(\Delta)
$$
where $\bar{u}'_{\textbf{T}}$ is the composition 
$$
\bar{\mathcal{C}}'_0 \ \xrightarrow{\ \bar{u}'\ } \ \bar{\mathcal{P}}' \times_{\textbf{T} \times \textbf{G}} \tilde{V} \ \rightarrow \ \bar{\mathcal{P}}' \times_{\textbf{T} \times \textbf{G}} W\ \cong\ \bar{\mathcal{Q}}' \times_{\textbf{T} } W.
$$
By adding additional marked points $\bar{v}_1,\cdots, \bar{v}_t: \Delta \rightarrow \bar{\mathcal{C}}'_0$ (it is possible since we have a family over $\Delta$, not $\Delta^0$) until 
$$
\left((\mathcal{C}, p_1,\cdots, p_k, v_1,\cdots, v_t),\ f: \mathcal{C} \rightarrow Y\right)\ \in\ \overline{\mathcal{M}}_{g,k+t}(Y, \beta')(\Delta^0)
$$
becomes stable
with $v_i=\bar{v}_i|_{\Delta^0}$, 
we have its extension over $\Delta$
\begin{align} \label{ex1}
\left((\bar{\bar{\mathcal{C}}}, \bar{\bar{p}}_1, \cdots, \bar{\bar{p}}_k, \bar{\bar{v}}_1,\cdots, \bar{\bar{v}}_t),\ \bar{\bar{f}}:\bar{\bar{\mathcal{C}}} \rightarrow Y\right) \ \in\  \overline{\mathcal{M}}_{g,k+t}(Y, \beta')(\Delta).
\end{align}
Let $\bar{\bar{\phi}}: (\bar{\bar{\mathcal{C}}}, \bar{\bar{p}}_1, \cdots, \bar{\bar{p}}_k) \rightarrow (\bar{\bar{\mathcal{C}}}_0, \bar{\bar{p}}_1, \cdots, \bar{\bar{p}}_k)$ be the contraction of rational tails.
Then it is obvious that 
$$
(\bar{\bar{\mathcal{C}}}_0, \bar{\bar{p}}_1, \cdots, \bar{\bar{p}}_k) |_{\Delta^0}\ \cong\ (\mathcal{C}_0, p_1, \cdots, p_k).
$$
Again, by adding marked points $\bar{\bar{r}}_1, \cdots, \bar{\bar{r}}_l$ on $\bar{\bar{\mathcal{C}}}_0$ until
$$
\left((\mathcal{C}_0,p_1,\cdots,p_k, r_1, \cdots, r_l),\ \mathcal{P},\ u\right)\ \in\ Q_{g,k+l}(\tilde{V} /\!\!/ \textbf{T} \times \textbf{G}, \beta'')(\Delta^0)$$
becomes stable with $r_i := \bar{\bar{r}}_i |_{\Delta^0}$, we can find its extension
\begin{align} \label{ex2}
\left((\tilde{\mathcal{C}}_0,\tilde{p}_1,\cdots,\tilde{p}_k, \tilde{r}_1, \cdots, \tilde{r}_l),\ \tilde{\mathcal{P}},\ \tilde{u}\right) \ \in\ Q_{g,k+l}(\tilde{V} /\!\!/ \textbf{T} \times \textbf{G}, \beta'')(\Delta).
\end{align}
Though \eqref{ex1} and \eqref{ex2} seem to give an extension by forgetting extra marked points, we need one more step for \eqref{ex1} since it may have redundant components.
For now, we cannot contract them yet.

To be more precise, adding marked points $\tilde{q}_1, \cdots, \tilde{q}_m$ on $\tilde{\mathcal{C}}_0$ until when
\begin{align}
\label{ext4}
\left((\mathcal{C},p_1,\cdots,p_k, q_1, \cdots, q_m),\ f \right)\ \in\ \overline{\mathcal{M}}_{g,k+m}(Y, \beta')(\Delta^0)
\end{align}
{\em and} 
\begin{align*}
& \left((\tilde{\mathcal{C}}_0,\tilde{p}_1,\cdots,\tilde{p}_k, \tilde{q}_1, \cdots, \tilde{q}_m),\ \tilde{\mathcal{Q}} := \tilde{\mathcal{P}}/\textbf{G}, \  \tilde{u}_{\textbf{T}}: \mathcal{C}_0 \rightarrow \mathcal{Q} \times_{\textbf{T}} W\right)\\
& \in \ Q_{g,k+m}(W /\!\!/_{\bf 1} \textbf{T} , d)(\Delta)
\end{align*}
become stable with $q_i := \tilde{q}_i |_{\Delta^0}$, we obtain an extension of \eqref{ext4}, 
$$
\left((\tilde{\mathcal{C}},\tilde{p}_1,\cdots,\tilde{p}_k, \tilde{q}_1, \cdots, \tilde{q}_m),\ \tilde{f}\right)\ \in\ \overline{\mathcal{M}}_{g,k+m}(Y, \beta')(\Delta).
$$ 
Now let $(\hat{\mathcal{C}}, \hat{p}_1, \cdots, \hat{p}_k)$ be {\em the contraction of rational components in rational tails} of $(\tilde{\mathcal{C}}, \tilde{p}_1, \cdots, \tilde{p}_k)$ on which $\tilde{f}$ is constant. Note that we could not do it before since we had the extension of $\cC'_0$, not of $\cC_0$.
We define $\hat{f}:\hat{\mathcal{C}} \rightarrow Y$ for which 
$$
\xymatrix@R=6mm{
\tilde{\mathcal{C}}\ar[r]\ar[rd]_{\tilde{f}} & \hat{\mathcal{C}} \ar[d]^{\hat{f}}\\
& Y
}
$$
is commutative. So, we have an extension of $((\mathcal{C}, p_1, \cdots, p_k),\ f)$,
\begin{align}
\label{ext3}
\left( (\hat{\mathcal{C}}, \hat{p}_1, \cdots, \hat{p}_k),\ \hat{f} \right)\ \in\ \mathfrak{M}'_{g,k}(Y, \beta')(\Delta).
\end{align}
Let $\hat{\phi}: (\hat{\mathcal{C}}, \hat{p}_1, \cdots, \hat{p}_k) \rightarrow (\hat{\mathcal{C}}_0, \hat{p}_1, \cdots, \hat{p}_k)$ be the contraction of rational tails.
Then it is obvious that 
$$
(\hat{\mathcal{C}}_0, \hat{p}_1, \cdots, \hat{p}_k) |_{\Delta^0} \ \cong \ (\mathcal{C}_0, p_1, \cdots, p_k).
$$
The construction gives rise to the contraction morphism
$$
c \ : \ (\hat{\mathcal{C}}_0, \hat{p}_1, \cdots, \hat{p}_k)\ \longrightarrow\ (\tilde{\mathcal{C}}_0, \tilde{p}_1, \cdots, \tilde{p}_k).
$$
So, we obtain
\begin{align}
\label{ggg}
\left( (\hat{\mathcal{C}}_0,\hat{p}_1,\cdots,\hat{p}_k,),\ \hat{\mathcal{P}} := c^*\tilde{\mathcal{P}},\ \hat{u} := c^*\tilde{u} \right)\ \in\ Q^{pre}_{g,k}(\tilde{V} /\!\!/ \textbf{T} \times \textbf{G}, \beta'')(\Delta).
\end{align}
We can check the pair (\eqref{ggg}, \eqref{ext3}) is in $Q^{pre}_{g,k}(E /\!\!/ \textbf{G}, \beta)(\Delta)$. 
We can define its stabilisation using the ample line bundle $\mathcal{L}$ \eqref{lb} on $E/\!\!/ \textbf{G}$.
This produces an extension of \eqref{data}.
\end{proof}

\subsection{Obstruction theory} \label{PoT}

For a reductive group $\bf{H}$, let $\mathfrak{B}\text{un}_{\textbf{H}}^{g,k}$ be the moduli space consisting of 
$((C_0,p_1, \cdots, p_k),P)$,
where $(C_0,p_1, \cdots, p_k)$ is a genus $g$, $k$-pointed prestable curve, and $P$ is a principal $\textbf{H}$-bundle on $C_0$.
Let $\mathfrak{S}_{g,k,\beta_0}$ be the stack consisting of $((C_0,p_1, \cdots, p_k),P,u \in \Gamma(C_0, P \times_{\textbf{T} \times \textbf{G}} V))$ where $((C_0,p_1, \cdots, p_k),P) \in \mathfrak{B}\text{un}_{\textbf{T}\times \textbf{G}}^{g,k}$ such that $u^{-1}(P \times_{\textbf{T} \times \textbf{G}} V^{un}(\textbf{G}, \theta) )$ is a set of finite points in $C_0^{sm} \smallsetminus \{p_1, \cdots, p_k \}$. 
Then we have a fiber diagram of forgetting morphisms
$$
\xymatrix@R=6mm{
Q^{pre}_{g,k}(\widetilde{V}/\!\!/_{\widetilde{\theta}} (\textbf{T} \times \textbf{G}) , \beta'') \ar[r]\ar[d]  &  \mathfrak{S}_{g,k,\beta_0} \ar[d]\\
Q^{pre}_{g,k}(W /\!\!/_{\bf 1} \textbf{T}, d) \ar[r]  &  \mathfrak{B}\text{un}_{\textbf{T}}^{g,k}
},
$$
and hence the following is a fiber diagram by \eqref{fiberProd},
\begin{align}
\label{prediag}
\xymatrix@R=6mm{
Q^{pre}_{g,k}(E /\!\!/ \textbf{G},\beta)\ar[r]\ar[d] & Q^{pre}_{g,k}(\widetilde{V} /\!\!/_{\widetilde{\theta}} (\textbf{T} \times \textbf{G}) , \beta'') \ar[d]\ar[r]   &\mathfrak{S}_{g,k,\beta_0} \ar[d]\\
\mathfrak{M}'_{g,k}(Y,\beta') \ar[r] & Q^{pre}_{g,k}(W /\!\!/_{\bf 1} \textbf{T}, d) \ar[r] & \mathfrak{B}\text{un}_{\textbf{T}}^{g,k}.
}
\end{align}

Let $\mu_1 : \mathfrak{S}_{g,k,\beta_0} \rightarrow \mathfrak{B}\text{un}_{\textbf{T}}^{g,k}$ be the right most vertical morphism in \eqref{prediag}. This can be factored into 
$$
\mathfrak{S}_{g,k,\beta_0} \ \xrightarrow{\ \mu_1^1\ }  \ \mathfrak{B}\text{un}_{\textbf{T} \times \textbf{G}}^{g,k} \ \xrightarrow{\ \mu_1^2 \ } \ \mathfrak{B}\text{un}_{\textbf{T}}^{g,k}.
$$
Consider the universal curve $\mathfrak{C}_0$, the universal principal $\textbf{T} \times \textbf{G}$-bundle $\mathfrak{P}$, and the universal section $u$, 
$$
\xymatrix@R=6mm{
\mathfrak{P} \times_{(\textbf{T} \times \textbf{G})} V \ar[d]^-{\rho} &  \\
\mathfrak{C}_0   \ar[d]^-{\pi_1} \ar@/^/[u]^-{u} & \\
\mathfrak{S}_{g,k,\beta_0}  \ar[r]^-{\mu_1^1}&  \mathfrak{B}\text{un}_{\textbf{T} \times \textbf{G}}^{g,k}.
}
$$
By using the similar argument of the proof of \cite[Theorem 4.5.2]{CKM}, we obtain a relative perfect obstruction theory for $\mu_1^1$
$$
E_{1}^{\bullet} \ :=\ \left(R^\bullet (\pi_1)_*(u^* T_\rho)\right)^\vee\ \in\ D^b (\mathfrak{S}_{g.k,\beta_0}),
$$
where $T_\rho$ is the relative tangent complex of $\rho$. 
We define a complex
$$
E_{\mu_1}^{\bullet}\ :=\ \text{Cone}\left(E^{\bullet}_1[-1] \rightarrow \mathbb{L}_{\mu_1^1}[-1] \rightarrow (\mu_1^1)^*\mathbb{L}_{\mu_1^2}\right)\ \in\ D^b (\mathfrak{S}_{g.k,\beta_0}),
$$
where $\mathbb{L}_{\mu_1^1}$ (resp. $\mathbb{L}_{\mu_1^2}$) is the relative cotangent complex for $\mu_1^1$ (resp. $\mu_1^2$).

Let $\mu_2 : \mathfrak{M}'_{g,k}(Y, \beta') \rightarrow \mathfrak{B}\text{un}_{\textbf{T}}^{g,k}$ be the composition of the bottom morphisms in \eqref{prediag}. 
Consider the composition
$$
\mu' \ :\  \mathfrak{M}'_{g,k}(Y, \beta') \ \xrightarrow{\ \mu_2 \ } \ \mathfrak{B}\text{un}_{\textbf{T}}^{g,k} \ \xrightarrow{\ \mu_3\ } \ \text{Spec}\mathbb{C}.
$$
It can be factored into 
$$
\mathfrak{M}'_{g,k}(Y, \beta') \ \xrightarrow{\ \mu'_1\ } \ \mathfrak{M}_{g,k}   \ \xrightarrow{\ \mu'_2\ }  \ \text{Spec}\mathbb{C},
$$
where $\mu'_1$ is the forgetting morphism. Consider the universal curve, and the universal morphism $f$,
$$
\xymatrix@R=6mm{
Y & \mathfrak{C}   \ar[d]^-{\pi_2} \ar[l]_-{f} & \\
& \mathfrak{M}'_{g,k}(Y,\beta')  \ar[r]^-{\mu'_1}&  \mathfrak{M}_{g,k}.
}
$$
The complex 
$$
E_{2}^{\bullet}\ :=\ \left(R^\bullet (\pi_2)_*(f^* T_Y)\right)^\vee \in D^b (\mathfrak{M}'_{g.k}(Y,\beta'))
$$ 
is then a relative perfect obstruction theory for $\mu'_1$.
We define a complex
$$
E'^{\bullet}\ :=\ \text{Cone}\left(E^{\bullet}_2[-1] \rightarrow \mathbb{L}_{\mu'_1}[-1] \rightarrow (\mu'_1)^*\mathbb{L}_{\mathfrak{M}_{g,k}}\right)\ \in\ D^b (\mathfrak{M}'_{g.k}(Y,\beta')).
$$
To find a perfect obstruction theory relative to $\mathfrak{B}\text{un}_{\textbf{T}}^{g,k}$ eventually, we would like to find a morphism $(\mu_2)^*\mathbb{L}_{\mathfrak{B}\text{un}_{\textbf{T}}^{g,k}} \ra E'^{\bullet}$ commuting the diagram
$$ 
\xymatrix@R=5mm{
(\mu_2)^*\mathbb{L}_{\mathfrak{B}\text{un}_{\textbf{T}}^{g,k}} \ar[r]\ar[rd] & E'^{\bullet} \ar[d]\\
& \mathbb{L}_{\mathfrak{M}'_{g,k}(Y,\beta')} .
}
$$ 
To do so, consider the contraction of rational tails $\phi: \fC \ra \fC_0$ of the universal curve on $\mathfrak{M}'_{g.k}(Y,\beta')$.
Then we obtain a commutative (not fiber) diagram
$$
\xymatrix@R=6mm{
T_{\fC} (-\sum_i p_i) \ar[r] \ar[d]_-{v_1} & f^*T_Y \ar[d]^{v_2} \\
\phi^*  T_{\fC_0} (-\sum_i p_i)  \ar[r] & \phi^* T_{B{\bf T}} |_{\fC_0}.
}
$$
Indeed, the complex $E'^\bullet$ is isomorphic to the dual derived push forward of 
$$
\mathrm{Cone}(v_1)[-1]\  \ra\ f^*T_Y .
$$
Hence the morphism is induced by $v_2$.
Now we define a complex
$$
E_{\mu_2}^{\bullet}\ :=\ \text{Cone}\left((\mu_2)^*\mathbb{L}_{\mathfrak{B}\text{un}_{\textbf{T}}^{g,k}} \rightarrow E'^{\bullet}\right)\ \in\ D^b (\mathfrak{M}'_{g.k}(Y,\beta')).
$$
 
Letting $p_1: Q^{pre}_{g,k}(E /\!\!/ \textbf{G},\beta) \rightarrow \mathfrak{S}_{g,k,\beta_0}$, $p_2 : Q^{pre}_{g,k}(E /\!\!/ \textbf{G},\beta) \rightarrow \mathfrak{M}'_{g,k}(Y,\beta')$ be the forgetful morphisms, we define 
\begin{align} \label{defthy}
E^\bullet \:= \ \left(p_1^* E_{\mu_1}^\bullet \oplus p_2^* E_{\mu_2}^\bullet \right) |_{Q_{g,k}(E /\!\!/ \textbf{G}, \beta)}\ \in\ D^b (Q_{g,k}(E /\!\!/ \textbf{G}, \beta)).
\end{align}
It is a relative perfect obstruction theory for $Q_{g,k}(E /\!\!/ \textbf{G}, \beta)$ over $\mathfrak{B}\text{un}_{\textbf{T}}^{g,k}$ \cite[Propositions 7.2, 7.4]{BF}.

\begin{Thm}
\label{constthm}
The moduli space $Q_{g,k}(E /\!\!/ \bf{G} ,\beta)$ is a proper DM stack over $\Spec \; \mathbb{C}.$ Moreover, it is equipped with a natural perfect obstruction theory.
\end{Thm}

The virtual fundamental class
$[Q_{g,k}(E /\!\!/ \textbf{G} ,\beta)]^{vir}$ has dimension
\begin{align} \label{vdim}
\beta(c_1([E /\textbf{G}])) + (1-g)(\text{dim}(E /\!\!/ \textbf{G})-3) +k.
\end{align}

\subsection{Graph spaces} \label{StableQGmap}

For $g,k \in \mathbb{Z}_{\geq 0}$ and $\beta=(\beta' , \beta_0)$, consider the following data:
\begin{align} \label{GraphQuasiMap}
\left(\phi:(C,p_1,\cdots,p_k) \rightarrow (C_0,p_1,\cdots,p_k), \ f:C \rightarrow Y,\ P, \ u,\ \alpha:C \rightarrow \mathbb{P}^1\right)
\end{align}
satisfying
\begin{itemize}
\item $(C,p_1, \cdots, p_k)$ and $(C_0,p_1, \cdots, p_k)$ are prestable $k$-pointed curves of genus $g$,
\item $\alpha$ is a regular map with degree $1 \in H_2(\mathbb{P}^1) \cong \mathbb{Z}$. This condition is equivalent to the fact that there is an irreducible rational component $C_1$ in $C$ (called {\em the distinguished component}) such that $\alpha|_{C_1}$ is an isomorphism and $\alpha$ on other components are contractions, 
\item $\phi$ is the contraction of all rational tails on $C$. By a rational tail here, we mean a maximal connected tree of rational curves not having $C_1$ as a component with no marked points, attached to other components at only one node on $C$,
\item $f$ is of degree $\beta'$ and $f$ restricted to each component of a rational tail is non-constant,
\item $P$ is a $\textbf{G}$-principal bundle on $C_0$, 
\item $u \in \Gamma\left(C_0, P_f \times_{(\textbf{T} \times \textbf{G})} V \right)$, 
\item $\beta_0 (L)= \text{deg}\left(u^*(P_f \times_{(\textbf{T} \times \textbf{G})} L)\right)$ for all $L \in \text{Pic}^{\textbf{T} \times \textbf{G}}(V)$.
\end{itemize}

\noindent We call the above data \eqref{GraphQuasiMap} {\em a genus $g$, $k$-pointed graph quasimap with degree $\beta$}, or {\em graph quasimap with type $(g,k,\beta)$}. 
It is {\em $\theta$-prestable} if $(C_0,P_f, \widetilde{u},\alpha_0)$ is a $\widetilde{\theta}$-prestable graph quasimap in the sense of \cite{CKM} where $\alpha_0:C_0 \rightarrow \mathbb{P}^1$ is the induced morphism by $\alpha$.  
It is {\em $\theta$-stable} if for each irreducible component $C' \subset C_0$ which is not $C_1$, $f|_{C'}$ is non-constant or $(C_0, P_f, \widetilde{u})|_{C'}$ is $0^+$-stable in the sense of \cite{CKM} with respect to $\widetilde{\theta}$.
We denote the moduli space of stable graph quasimaps with type $(g,k,\beta)$ by 
$QG_{g,k}(E /\!\!/ \textbf{G}, \beta)$. 

\begin{Thm}
The moduli space $QG_{g,k}(E /\!\!/ \bf{G} ,\beta)$ is a proper DM stack over $\Spec\; \mathbb{C}.$ Moreover, it is equipped with a natural perfect obstruction theory.
\end{Thm}

We define a $\mathbb{C}^*$-action on $QG_{g,k}(E /\!\!/ \bf{G} ,\beta)$
\begin{eqnarray*}
& \mathbb{C}^* \times QG_{g,k}(E /\!\!/ \textbf{G} ,\beta) \ \longrightarrow \ QG_{g,k}(E /\!\!/ \textbf{G} , \beta), \\
& (t, \ (\phi, f, P, u, \alpha)) \ \longmapsto \ (\phi, f, P, u, \alpha_t)
\end{eqnarray*}
where $\alpha_t$ is the composition of $\alpha$ with 
$\mathbb{P}^1 \rightarrow \mathbb{P}^1, ~~~[x;y] \mapsto [tx;y]$.

\subsection{$g=0$ graph spaces} \label{g0:graph}
For $m,k \in \mathbb{Z}_{\geq 0}$ and $\beta=(\beta' , \beta_0)$, consider the following data,
{\medmuskip=0mu
\thinmuskip=0mu
\thickmuskip=0mu
$$ 
\left(\phi\ :\ (C,q_1, \cdots, q_m, p_1,\cdots,p_k)\ \rightarrow\ (C_0,q_1,\cdots,q_m), \ f\ :\ C\ \rightarrow\ Y,\ P, \ u,\ \alpha\ :\ C \ \rightarrow\ \mathbb{P}^1\right)
$$
}
\par \noindent
satisfying
\begin{itemize}
\item $(C,q_1, \cdots, q_m, p_1, \cdots, p_k)$ is a $g=0$ prestable $(m+k)$-pointed curve,
\item $(C_0,q_1, \cdots, q_m)$ is a $g=0$ prestable $m$-pointed curve,
\item $\alpha$ is a regular map with degree $1 \in H_2(\mathbb{P}^1) \cong \mathbb{Z}$, 
\item $\phi$ is the contraction of all rational tails on $C$ \emph{after forgetting $p_1, \cdots, p_k$}, 
\item $f$ is of degree $\beta'$ and stable on each contracted component,
\item $P$ is a $\textbf{G}$-principal bundle on $C_0$, 
\item $u \in \Gamma \left(C_0, P_f \times_{(\textbf{T} \times \textbf{G})} V \right)$, 
\item $\beta_0 (L)= \text{deg}\left(u^*(P_f \times_{(\textbf{T} \times \textbf{G})} L)\right)$ for all $L \in \text{Pic}^{\textbf{T} \times \textbf{G}}(V)$.
\end{itemize}
We can define the stability condition for the above data as before to construct a moduli space denoted by $QG_{m|k}(E /\!\!/ \bf{G} , \beta)$.
Then it is a proper DM stack over $\Spec\; \CC$ equipped with a natural perfect obstruction theory.
Furthermore, it comes with a natural $\CC^*$-action.

There are evaluation maps $\text{ev}_a: QG_{m|k}(E /\!\!/ {\bf G} , \beta) \ra Y$ evaluating at $p_1, \cdots, p_k$, and $\text{ev}_b: QG_{m|k}(E /\!\!/ {\bf G} , \beta) \ra E /\!\!/ {\bf G}$ evaluating at $q_1, \cdots, q_m$.

\section{Characterisation properties for $I$-function} \label{PolyI}

In this section, we will show that $I^{{\bf S}, \mu}_{\cF^{(n)}} = \mu^*I^{\bf S}_{\cF^{(n)}}$ does not have poles outside of $z=0, \infty, -\chi_{\mu,\nu}/k$ and satisfies the initial and polynomiality conditions using quasimap moduli spaces.

\subsection{Graph space expression of $I^{\bf S}_{\cF^{(n)}}$ and its initial condition} \label{QuasiI}

Let 
$QF_{k,\beta} \subset QG_{0|k}({\cF^{(n)}} , \beta)^{\mathbb{C}^*}$
be a component of the fixed locus, where each node, marked point and degree is concentrated on $0 \in \mathbb{P}^1 \cong C_0$.
Since $\mathbb{P}^1 \backslash \{0\}$ maps constantly to ${\cF^{(n)}}$, we obtain an evaluation map
$\text{ev}_\bullet: QF_{k,\beta} \ra {\cF^{(n)}}$.

\begin{Prop} \label{Ishape}
The $I$-function $I^{\bf{S}}_{\cF^{(n)}}$ can be written as 
\begin{align*}
I^{\bf{S}}_{\cF^{(n)}} \ =\ e^{\frac{\bf{t}}{z}}  \sum_{k, \beta =(d,D)} \frac{q^d Q^D}{k!} e^{\sum_id_it_i}  (\mathrm{ev}_{\bullet})_* \left( \frac{[QF_{k,\beta}]^{vir} \prod_{a=1}^k \mathrm{ev}^*_a(u)}{e_{\mathbb{C}^* \times {\bf{S}}}(N^{vir}_{QF_{k,\beta}/ QG_{0|k}({\cF^{(n)}},\beta)})} \right). 
\end{align*}
\end{Prop}

\begin{proof}
Let $\widetilde{F}_{k,\beta} := F_{k,D} \times _{\overline{\cM}G_{0,k}(Y, D)   }   QG_{0|k}({\cF^{(n)}},\beta)$ be the fiber product sitting in the diagram
$$
\xymatrix@R=5mm{
QF_{k,\beta} \ar[r] & \widetilde{F}_{k,\beta} \ar[r] \ar[d] & QG_{0|k}({\cF^{(n)}},\beta) \ar[d] \\
& F_{k,D} \ar[r] & \overline{\cM}G_{0,k}(Y, D)
}
$$
where $F_{k,D} \hookrightarrow \overline{\cM}G_{0,k}(Y, D)$ is defined in \eqref{graphJform}. 
Using \eqref{defthy}, one can decompose the virtual normal bundle (with respect to the $\CC^*$-action) to get
{\medmuskip=0mu
\thinmuskip=0mu
\thickmuskip=0mu
\begin{align} \label{qeqe}
e_{\mathbb{C}^* \times \bf{S}}(N^{vir}_{QF_{k,\beta}/ QG_{0|k}({\cF^{(n)}},\beta)})\ =\ e_{\mathbb{C}^* \times \bf{S}}( N^{vir}_{F_{k,D} /\overline{\mathcal{M}}G_{0,k}(Y  , D)})\cdot e_{\mathbb{C}^* \times \bf{S}}(N^{vir}_{QF_{k,\beta} / \widetilde{F}_{k,\beta}}).
\end{align}
}
\par \noindent
Note that the $\bf{S}$-action on $N^{vir}_{F_{k,D} /\overline{\mathcal{M}}G_{0,k}(Y  , D)}$ is trivial.
As explained in \eqref{graphJform}, we have
{\medmuskip=0mu
\thinmuskip=0mu
\thickmuskip=0mu
$$e_{\mathbb{C}^* \times \bf{S}}( N^{vir}_{F_{k,D} /\overline{\mathcal{M}}G_{0,k}(Y  , D)})\ =\ e_{\mathbb{C}^*}( N^{vir}_{F_{k,D} /\overline{\mathcal{M}}G_{0,k}(Y  , D)})\ =\
\left\{
\begin{array}{cl}
1 & \text{if } \beta' =0 \\
z(z-\psi) & \text{if } \beta' \neq 0
\end{array}
\right.$$ 
}
\par \noindent
which contributes to $\pi^*(J_{D})$.
Since the flag bundle $\cF^{(n)}$ can be thought of as a tower of Grassmannian bundles, $QF_{k, \beta}$ is a disjoint union of products of flag bundles on $Y$ \cite[Lemma 1.2]{BCK1}.
Following the computation in the proofs of \cite[Theorem 1.5]{BCK1} and \cite[Theorem 1]{BCK2}, we obtain
\begin{align*} 
e_{\mathbb{C}^* \times \bf{S}}(N^{vir}_{QF_{k,\beta} / \widetilde{F}_{k,\beta}})^{-1} & = \sum_{\sum_l d_i^l =d_i} \prod_{i=1}^n  \left( \prod_{1 \leq l \neq l' \leq r_i}\frac{\prod_{s=-\infty}^{d_i^l -d_i^{l'}}(H_{i,l}-H_{i,l'}+sz) }{\prod_{s=-\infty}^{0}(H_{i,l}-H_{i,l'}+sz)}  \right. \nonumber \\
& \left.  \prod_{1 \leq l \leq r_i,~ 1 \leq l' \leq r_{i+1}}\frac{\prod_{s=-\infty}^{0}(H_{i,l}-H_{i+1,l'}+sz)}{\prod_{s=-\infty}^{d_i^l -d_{i+1}^{l'}}(H_{i,l}-H_{i+1,l'}+sz)} \right).
\end{align*}
More precisely, $\mathrm{ev}^*_\bullet$ of the right-hand side is equal to the left-hand side.
\end{proof}

As a corollary of Proposition \ref{Ishape}, we obtain the initial condition for $I^{\bf S}_{\cF^{(n)}}$.

\begin{Cor}
$I^{\bf S}_{\cF^{(n)}}$ satisfies the initial condition.
\end{Cor}
\begin{proof}
Let $\iota_\mu: QF_{k,\beta}^\mu \hookrightarrow QF_{k,\beta}$ be the substack taking the marked point to $Y^\mu \subset \cF^{(n)}$.
Then the pullback of the $I$-function $I^{{\bf S}, \mu}_{\cF^{(n)}} = \mu^*I^{\bf S}_{\cF^{(n)}}$ is
\begin{align*}
I^{\bf{S},\mu}_{\cF^{(n)}} \ =\ e^{\frac{\mu^*\bf{t}}{z}}  \sum_{k, \beta =(d,D)} \frac{q^d Q^D}{k!} e^{\sum_id_it_i}  (\mathrm{ev}_{\bullet})_* \left( \frac{[QF^{\mu}_{k,\beta}]^{vir} \prod_{a=1}^k \mathrm{ev}^*_a(u)}{e_{\mathbb{C}^* \times {\bf{S}}}(N^{vir}_{QF^{\mu}_{k,\beta}/ QG_{0|k}({\cF^{(n)}},\beta)})} \right). 
\end{align*}
Since $N_{QF^{\mu}_{k,\beta}/QF_{k,\beta}} = ev_\bullet^*N_{Y^\mu/\cF^{(n)}}$, where $ev_\bullet$ denotes the evaluation map to $Y \cong Y^\mu$ (beware that $\mathrm{ev}_\bullet$ is the map to $\cF^{(n)}$), $I^{\bf{S},\mu}_{\cF^{(n)}}$ becomes
\begin{align*}
I^{\bf{S},\mu}_{\cF^{(n)}} \ =\ e^{\frac{\mu^*\bf{t}}{z}}  \sum_{k, \beta =(d,D)} \frac{q^d Q^D}{k!} e^{\sum_id_it_i}  (ev_{\bullet})_* \left( \frac{[QF^{\mu}_{k,\beta}]^{vir} \prod_{a=1}^k \mathrm{ev}^*_a(u)}{\iota^*_\mu e_{\mathbb{C}^* \times {\bf{S}}}(N^{vir}_{QF_{k,\beta}/ QG_{0|k}({\cF^{(n)}},\beta)})} \right). 
\end{align*}
Hence by using the decomposition \eqref{qeqe} as well as the string and divisor equations of Gromov--Witten theory, we observe that $(I^{\bf S}_{\cF^{(n)}})^\mu_Y$ (defined in Section \S \ref{subsub:initial}) becomes
\begin{align} \label{IYmu}
(I^{\bf S}_{\cF^{(n)}})^\mu_Y = \sum_{D} Q^D J_{D}(z, \mu^*{\bf t} + u) \cdot \iota_\mu^* \ e_{\mathbb{C}^* \times \bf{S}}\left(N^{vir}_{QF_{k,(d_{D,\mu},D)} / \widetilde{F}_{k,(d_{D,\mu},D)}}\right)^{-1}.
\end{align}

Recall that $\tilde{V}$, $W$, $\bf G$ and $\bf T$ are introduced in Section \S \ref{Sect:pr} for the GIT presentation $E /\!\!/ {\bf G}$ of $\cF^{(n)}$.
Letting $\pi$ and $f$ be the universal ones
$$
\xymatrix@R=5mm{
QF_{k,\beta} \times \PP^1 \ar[r]^f \ar[d]_\pi & [\tilde{V}/{\bf T} \times {\bf G}]\\
QF_{k,\beta} ,
}
$$
we obtain an isomorphism 
\begin{align} \label{NTaaa}
N^{vir}_{QF_{k,\beta} / \widetilde{F}_{k,\beta}} \ \cong \ (R\pi_*f^*T_{[\tilde{V}/{\bf T} \times {\bf G}]/[W/{\bf T}]})^{\mathrm{mov}}
\end{align}
by the construction of $E^\bullet_{\mu_1}$ from \eqref{prediag} in Section \S \ref{PoT}.
Using the notations in \eqref{V}, we define a ${\bf T} \times {\bf G}$-equivariant subspace of $V$
$$
V^\mu\ :=\ \prod_i GL\left(\prod_{j \in I^\mu_i} \CC_{e_j}\right) \ \subset\ V.
$$
Note that the $\bf G$-action on $V^\mu$ is given by the right multiplications.
Let $[W/{\bf T}]^\mu := [W \times V^{\mu}/{\bf T} \times {\bf G}] \subset [\tilde{V}/{\bf T} \times {\bf G}]$ be the corresponding substack. Then it is the $\bf S$-fixed locus in $[\tilde{V}/{\bf T} \times {\bf G}]$ corresponding to $\mu$, isomorphic to $[W/{\bf T}]$.
For a morphism $\PP^1 \ra [W/{\bf T}]$ of degree $D$, the pullback of the line bundle $\det[W \times \left(\prod_{j \in I^\mu_i} \CC_{e_j}\right) / {\bf T}]$ to $\PP^1$ is of degree $(d_{D,\mu})_i$.
Hence when $d=d_{D,\mu}$, the universal map $f$ factors through 
$$
f: QF_{k,(d_{D,\mu},D)}^\mu \times \PP^1 \ \longrightarrow\ [W/{\bf T}]^\mu  \ \hookrightarrow \ [\tilde{V}/{\bf T} \times {\bf G}].
$$
Since $N_{[W/{\bf T}]^\mu/[\tilde{V}/{\bf T}\times{\bf G}]} \cong T_{[\tilde{V}/{\bf T} \times {\bf G}]/[W/{\bf T}]} |_{[W/{\bf T}]^\mu}$, we have
$$
\iota_\mu^* N^{vir}_{QF_{k,(d_{D,\mu},D)} / \widetilde{F}_{k,(d_{D,\mu},D)}} \ \cong \ (R\pi_*f^*N_{[W/{\bf T}]^\mu/[\tilde{V}/{\bf T}\times{\bf G}]})^{\mathrm{mov}}
$$
by \eqref{NTaaa}.
Its Euler class is then
\begin{align} \label{inverseE}
& \iota_\mu^* \  e_{\mathbb{C}^* \times \bf{S}}\left(N^{vir}_{QF_{k,(d_{D,\mu},D)} / \widetilde{F}_{k,(d_{D,\mu},D)}}\right)^{-1}  \\ 
& \qquad \qquad = \ e_{\mathbb{C}^* \times \bf{S}}\left((R\pi_*f^*N_{[W/{\bf T}]^\mu/[\tilde{V}/{\bf T}\times{\bf G}]})^{\mathrm{mov}}\right)^{-1}  \nonumber \\
& \qquad \qquad = \ \prod_{i=1}^n
\prod_{ l \in I_i^\mu ,~  l' \in I_{i+1}^\mu}\frac{\prod_{s=-\infty}^{0}(H_{l}-H_{l'}+sz)}{\prod_{s=-\infty}^{d_l -d_{l'}}(H_{l}-H_{l'}+sz)}, \nonumber
\end{align}
where $H_l = - c_1(L_l)$ and $d_l = D \cap H_l$. 
Putting \eqref{inverseE} to \eqref{IYmu}, $(I^{\bf S}_{\cF^{(n)}})^\mu_Y(-z)$ lies on $-z^{-1}\cL^\mu_Y$ by \cite[Theorem $2^\prime$]{CG}.
\end{proof}

\subsection{Recursion relation for $I^{\bf S}_{\cF^{(n)}}$}

There is another description of $I^{\bf S}_{\cF^{(n)}}$ following the idea of \cite[Theorem 5.4.1]{CK0}.
Let 
$QF_{1|k,\beta} \subset QG_{1|k}({\cF^{(n)}} , \beta)^{\mathbb{C}^*}$
be a component of the fixed locus, where each node, marked point and degree is concentrated on $0 \in \mathbb{P}^1$.
Let $p_0:= e_{\CC^*}(\cO(1) \ot \CC_1)$, $p_\infty := e_{\CC^*}(\cO(1)) \in H^*_{\CC^*}(\PP^1)$, hence they satisfy 
$$
p_0 |_0 \ =\ z \ =\ - p_\infty|_\infty\ \  \text{ and } \ \ p_0 |_\infty \ =\ 0 \ =\  p_\infty|_0.$$
Let $S_{quasi}^*$ be the operator on $H^*_{\bf S}(\cF^{(n)},\QQ)(z) \ot_\QQ \QQ[\![q, Q]\!] $ 
{\medmuskip=-2mu
\thinmuskip=-2mu
\thickmuskip=-2mu
\begin{align} \label{Squasi}
S_{quasi}^*(\gamma) &\ :=\ e^{\frac{\bf t}{z}}\sum_{k, \beta =(d,D) } \frac{q^d Q^D}{k!}   e^{\sum_i d_i t_i}  (\text{ev}_{\bullet})_* \left( \frac{[QF_{1|k,\beta}]^{vir} \text{ev}_b^*(\gamma p_0) \prod_{a=1}^k \text{ev}^*_a(u)}{e_{\mathbb{C}^* \times {\bf S} }(N^{vir}_{QF_{1|k,\beta}/ QG_{1|k}({\cF^{(n)}},\beta)})} \right)\ ,
\end{align}
}
\par \noindent
and let $P_{quasi}$ be the series in $H^*_{\bf S}(\cF^{(n)},\QQ)[z] \ot_\QQ \QQ[\![q, Q, {\bf t}, u]\!] $ 
{\medmuskip=0mu
\thinmuskip=0mu
\thickmuskip=0mu
\begin{align} \label{Pquasi}
P_{quasi} & \ :=\ \sum_i \gamma_i \sum_{k,\beta =(d,D)} \frac{q^d Q^D}{k!}  e^{\sum_i d_i t_i}  \int_{[QG_{1|k}({\cF^{(n)}} , \beta)]^{vir}}  \text{ev}_b^*(\gamma^i p_\infty) \prod_{a=1}^k \text{ev}^*_a(u) . 
\end{align}
}

\begin{Prop} \label{Bfact}
The $I$-function $I^{\bf{S}}_{\cF^{(n)}}$ can be written as 
\begin{align*}
I^{\bf S}_{\cF^{(n)}} \ =\ S^*_{quasi} (P_{quasi}).
\end{align*}
\end{Prop}

\begin{proof}
The proof is identical to the proof of \cite[Theorem 5.4.1]{CK0}.
We include it here because it is omitted there.
By \cite[Proposition 5.3.1]{CK0}, one can check that the inverse of $S_{quasi}^*$ is
{\medmuskip=0mu
\thinmuskip=0mu
\thickmuskip=0mu
\begin{align*}
\gamma \ \mapsto\ & e^{-\frac{\bf t}{z}}\sum_i \gamma _i \sum_{k, \beta =(d,D)} \frac{q^d Q^D}{k!}  e^{\sum_i d_it_i}  \int_{[QF_{1|k,\beta}]^{vir}} \frac{\text{ev}_b^*(\gamma^i p_\infty) \text{ev}_{\bullet}^*(\gamma)\prod_{a=1}^k \text{ev}^*_a(u)}{e_{\mathbb{C}^* \times {\bf S} }(N^{'vir}_{1|k, \beta})}.
\end{align*}
}
\par \noindent
On the other hand, using virtual $\CC^*$-localisation \cite{GP}, one can check that 
\begin{align*}
& \int_{[QG_{1|k}(\cF^{(n)}, \beta)]^{vir}} \text{ev}_b^*(\gamma^i p_\infty) \prod_{a=1}^{k} \text{ev}^*_a ( u ) \\
=&  \sum_{\tiny{\begin{array}{c}\tiny{k_1 +k_2 =k}\\ \tiny{\beta_1 +\beta_2=\beta}\end{array}}} \int_{[QF_{k_1 ,\beta_1} \times_{\cF^{(n)}} QF_{1|k_2, \beta_2}]^{vir}} \frac{\text{ev}_b^*(\gamma^i p_\infty) \prod_{a=1}^{k} \text{ev}^*_a (u)}{e_{\mathbb{C}^* \times {\bf S}}(N_{k_1, \beta_1}^{vir})e_{\mathbb{C}^* \times {\bf S}}(N_{1|k_2, \beta_2}^{'vir})}\\
=&  \sum_{\tiny{\begin{array}{c}\tiny{k_1 +k_2 =k}\\ \tiny{\beta_1 +\beta_2=\beta}\end{array}}} \binom{k}{k_1 } \int_{[QF_{k_1 ,\beta_1}]^{vir}}\frac{\text{ev}^*_{\bullet}(\gamma^l) \prod_{a=1}^{ k_1} \text{ev}^*_a (u)}{e_{\mathbb{C}^* \times {\bf S}}(N_{k_1, \beta_1}^{vir})} \\
& \qquad\qquad\qquad\qquad\qquad \times  \int_{[QF_{1|k_2  ,\beta_2}]^{vir}}\frac{\text{ev}_{b}^*(\gamma^i p_\infty)\text{ev}^*_{\bullet}(\gamma_l) \prod_{a=1}^{k_2 } \text{ev}^*_a (u)}{e_{\mathbb{C}^* \times {\bf S}}(N_{1|k_2 , \beta_2}^{'vir})},\\
\end{align*}
where $N^{vir}_{m|k,\beta}=N^{vir}_{QF_{m|k, \beta}/QG_{m|k}(\cF^{(n)}, \beta)}$ and $N'$ denotes the $\CC^*$-bundle $N$ with the inverse $\CC^*$-action.
Hence by Proposition \ref{Ishape} we obtain $P_{quasi} = (S_{quasi}^*)^{-1}(I^{\bf S}_{\cF^{(n)}})$, and it induces
$$
I^{\bf S}_{\cF^{(n)}} \ =\  S^*_{quasi} (P_{quasi}).
$$ 
\end{proof}

Together with Proposition \ref{Prop:recI1}, the following corollary of Proposition \ref{Bfact} shows that $I^{\bf S}_{\cF^{(n)}}$ satisfies the recursion relation. 

\begin{Cor} \label{Cor:recI2}
$I^{{\bf S}, \mu}_{\cF^{(n)}}$ has poles only at $z=0, \infty$, and $-\chi_{\mu,\nu}/k$.
\end{Cor}
\begin{proof}
Since the denominators in $S^*_{quasi}(P_{quasi})$ are
\begin{align} \label{virnormal}
e_{\mathbb{C}^* \times {\bf S}}(N^{vir}_{1|k,\beta}) \ =\  \left\{
\begin{array}{cl}
z & \text{if }  \beta=0   \\
z(z-\psi) & \text{if } \beta \neq 0 ,
\end{array}
\right.
\end{align}
poles other than $z=0$, $\infty$ are $z=\psi$.
On the $\bf S$-fixed locus contributing to $I^{{\bf S}, \mu}_{\cF^{(n)}}$, 
$\psi$ can be either nilpotent (when the domain component containing marked point $\bullet$ lies on $Y^\mu$), or the pullback of $-\chi_{\mu,\nu}/k$ (when the domain component containing marked point $\bullet$ maps to the fiber). 
\end{proof}

\subsection{Polynomiality condition for $I^{\bf S}_{\cF^{(n)}}$}

The following proposition shows the polynomiality condition for $I^{\bf S}_{\cF^{(n)}}$.

\begin{Prop}
$I^{\bf S}_{\cF^{(n)}}$ satisfies the polynomiality condition.
\end{Prop}

\begin{proof}
The proof is parallel to that of Proposition \ref{Prop:polyJ}. Let 
$$
QF_{k_1, \beta_1}^{k_2, \beta_2}\ \subset\ QG_{0|k+1}({\cF^{(n)}}, \beta)^{\CC^*}
$$ 
be a component of the $\CC^*$-fixed locus, where $k_1$ marked points, degree $\beta_1$ are concentrated on $0 \in \mathbb{P}^1$ and $k_2$ marked points, degree $\beta_2$ are concentrated on $\infty \in \mathbb{P}^1$. Let 
$$
QG_{0|k+1}({\cF^{(n)}}, \beta)_{\mu} \subset QG_{0|k+1}({\cF^{(n)}}, \beta)^{\bf S}
$$ 
be the $\bf S$-fixed locus where the image of $\mathbb{P}^1$ lies on $Y^\mu$. 

We can construct a $\mathbb{C}^* \times \bf{S}$-equivariant line bundle $Q\cE_{i,\beta,\mu}$ on $QG_{0|k}({\cF^{(n)}}, \beta)$ for each $(i,\beta,\mu)$ such that
$$
Q\cE_{i,\beta,\mu} |_{QF_{k_1,\beta_1}^{k_2, \beta_2} \cap QG_{0|k+1}({\cF^{(n)}}, \beta)_{\mu} } \ =\ \mathbb{C}_{\beta_2 (\det \cE_i^\vee) +\pi_*\beta_2( \mu^*\det \cF_i )}.
$$
Letting $QN^{vir}_\mu := N^{vir}_{QG_{0|k+1}({\cF^{(n)}}, \beta)_{\mu} / QG_{0|k+1}({\cF^{(n)}}, \beta)}$, it follows from virtual $\CC^*$-lcoalisation \cite{GP} that the push forward of
\begin{align*}
 QZ_{\mu,j} := \displaystyle \sum_{k, \beta=(d,D) } \frac{q^d Q^D}{k!}  e^{\sum d_it_i} &  \frac{[QG_{0|k+1}({\cF^{(n)}}, \beta)_{\mu}]^{vir}}{e_{\CC^* \times \textbf{S}}(QN^{vir}_\mu)} \\
  &  \cap e^{\sum_i c_1(Q\cE_{i,\beta,\mu})y_i} \prod_{a=1}^k \text{ev}_a^*(u) \text{ev}^*_{k+1}(p_0 \delta_j ). 
\end{align*}
by $\mathrm{pt} : QG_{0|k+1}({\cF^{(n)}}, \beta)_{\mu} \to \Spec\; \CC$ becomes
\begin{align*}
\left( z\partial_{u_j} I^\mu (z,q),\  I^\mu (-z,qe^{-z\sum_i y_i E_i} )\right)_Y  \ =\ 
(\text{pt})_*QZ_{\mu,j}.
\end{align*}
Note that $QZ_{\mu,j}$ has no poles in $z=0$. 
This proves the polynomiality condition for $I^{\bf S}_{\cF^{(n)}}$.
\end{proof}

\section{Application} \label{Sect:App}

One can naturally ask if the Gromov--Witten theory of $\cF^{(n)}$ is related to that of $Y$.
There is one simple situation in which they are related each other -- $g=0$ invariants for semi-positive $\cF^{(n)}$, that is, 
$$ 
c_1(T_{\cF^{(n)}}) \ \cap\ \beta\ \geq\ 0 
$$
for any effective class $\beta$. 
We will discuss it in this section.
Throughout the section, we will forget $\bf S$-actions everywhere. 

\smallskip

The theorem below follows the idea in \cite[Section 5.5]{CK0}.
\begin{Thm}
Assume that $\cF^{(n)}$ is semi-positive. Then we obtain
\begin{align} \label{MirrorThm}
 I_{\cF^{(n)}} \ =\ I_0(qe^{\bf t},Q) \cdot J_{\cF^{(n)}} |_{{\bf t} + \pi^* u \mapsto \frac{{\bf t} + \pi^* u + I_1(qe^{\bf t},Q)}{I_0(qe^{\bf t},Q)}}
\end{align}
where $I_0(q,Q)$ is the coefficient of $1 \in H^0(\cF^{(n)}, \QQ[\![q,Q]\!])$ in $P_{quasi}|_{{\bf t}=u=0}$ and 
\begin{align} \label{I1exp}
I_1(q,Q)  \ := \ I_0(q,Q) \cdot \sum_{\beta = (d,D) \neq 0} q^dQ^D (\mathrm{ev}_1)_* \left( [Q_{0,2}(\cF^{(n)}, \beta)]^{vir} \mathrm{ev}^*_2(1) \right).
\end{align}
Here, $I_0(q,Q)$ is independent of $z$, invertible and 
$$
I_1(q,Q)\ \in\ H^{\leq 2} \left(\cF^{(n)}, \QQ[\![q,Q]\!] \right)
$$ 
so that the transformation
$$
{\bf t} + \pi^* u \ \longmapsto\ \frac{{\bf t} + \pi^* u + I_1(qe^{\bf t},Q)}{I_0(qe^{\bf t},Q)}
$$ 
makes sense.
Moreover, if $\cF^{(n)}$ is Fano of index $1$, namely, $c_1(T_{\cF^{(n)}}) \cap \beta \geq 1 $ for any effective class $\beta \neq 0$, then $I_0(q,Q) =1$.
If $\cF^{(n)}$ is Fano of index $2$, then $I_0(q,Q)=1$ and $I_1(q,Q)=0$ so that \eqref{MirrorThm} becomes $I_{\cF^{(n)}} = J_{\cF^{(n)}}.$
\end{Thm}

\begin{proof}
First of all, the virtual dimension of $QG_{1|0}(\cF^{(n)}, \beta)$ \eqref{vdim} is 
$$\beta(T_{\cF^{(n)}}) + \dim \cF^{(n)} +1 $$
which is greater or equal to $ \dim \cF^{(n)} +1$ by the semi-positivity condition. 
Since the maximal degree of $\gamma^i p_\infty$ in \eqref{Pquasi} is $\dim \cF^{(n)} +1,$ we obtain from \eqref{Pquasi} that
\begin{align} \label{Pexp} 
P_{quasi}|_{{\bf t}=u=0} = \left( \sum_{\tiny{\begin{array}{c}\tiny{\beta=(d,D)}\\ \tiny{\beta(T_{\cF^{(n)}}) =0 }\end{array}} } q^d Q^D  \int_{[QG_{1|0}({\cF^{(n)}} , \beta)]^{vir}}  \text{ev}_b^*\left(PD(pt) p_\infty \right) \right) \cdot 1.
\end{align}
Again by degree reasons, it does not have positive $z$-terms.
Hence we have
$P_{quasi}|_{{\bf t}=u=0} = I_0(q,Q) \cdot 1.$
Note that $I_0(q,Q) =  1 + O(q,Q)$, hence it is invertible.

By \eqref{Squasi}, \eqref{virnormal}, we obtain an asymptotic property of $S^*_{quasi}|_{{\bf t}=u=0}$,
{\medmuskip=-2mu
\thinmuskip=-2mu
\thickmuskip=-2mu
\begin{align} \label{aaaa}
S^*_{qausi}(\gamma)\ =\ \gamma\ +\ \frac{1}{z} \sum_{\beta = (d,D) \neq 0} q^dQ^D (\text{ev}_1)_* \left( [Q_{0,2}(\cF^{(n)}, \ \beta)]^{vir} \text{ev}^*_2(\gamma) \right) \ +\ O\left(\frac{1}{z^2}\right).
\end{align}}
\par \noindent
By applying $\gamma = P_{quasi}|_{{\bf t}=u=0}$ to \eqref{aaaa}, we obtain
\begin{align} \label{Iexp}
I_{\cF^{(n)}}|_{{\bf t}=u=0} \ =\ I_0(q,Q) \cdot 1\; +\; \frac{I_1(q,Q)}{z} \; +\; O\left(\frac{1}{z^2} \right)
\end{align}
by Proposition \ref{Bfact}.
Since the virtual dimension of $Q_{0,2}(\cF^{(n)}, \beta)$ \eqref{vdim} is 
\begin{align} \label{Q02dim}
\beta(T_{\cF^{(n)}}) \ +\ \dim \cF^{(n)}\ -\ 1 ,
\end{align}
we have $I_1(q,Q) \in  H^{\leq 2} \left(\cF^{(n)}, \QQ[\![q,Q]\!] \right)$.

Since $I_{\cF^{(n)}}|_{u=0} = e^{\frac{\bf t}{z}}I_{\cF^{(n)}}|_{{\bf t}=u=0, q \mapsto qe^{\bf t}} $, we obtain
\begin{align*}
I_{\cF^{(n)}} \ =\ I_0(qe^{\bf t},Q) \cdot 1 \; +\; \frac{1}{z} \left({\bf t} + \pi^* u \; +\; I_1(qe^{\bf t},Q)\right) +   O\left(\frac{1}{z^2}\right)
\end{align*}
from \eqref{Iexp}.
Then \eqref{MirrorThm} follows from the fact that $I_{\cF^{(n)}}$ lies on $-z^{-1}\cL ag_{\cF^{(n)}}$ (Theorem \ref{Mthm}) which is spanned by $J_{\cF^{(n)}}$.

When $\cF^{(n)}$ is Fano of index $1$, we have $I_0(q,Q)=1$ by observing the range of the summation in \eqref{Pexp}. When $\cF^{(n)}$ is Fano of index $2$, we have $I_1(q,Q)=0$ by the degree counting of \eqref{I1exp} since $\eqref{Q02dim} > \dim \cF^{(n)} $.
\end{proof}

\end{document}